\documentclass[11pt]{article}
 
\usepackage{amssymb,epsfig} 
\usepackage{amsmath}

\usepackage[T1]{fontenc}
\usepackage[dvips]{color}
\usepackage{amssymb}

\newtheorem{defi}{Definition}[section]
\newtheorem{prop}[defi]{Proposition}
\newtheorem{theo}[defi]{Theorem}
\newtheorem{conj}[defi]{Conjecture}
\newtheorem{lemm}[defi]{Lemma}
\newtheorem{coro}[defi]{Corollary}
\newtheorem{rema}[defi]{Remark}
\newtheorem{exem}[defi]{Example}
\newtheorem{exems}[defi]{Examples}

\newcommand{\bdefi}{\begin{defi}}
\newcommand{\edefi}{\end{defi}}
\newcommand{\bprop}{\begin{prop}}
\newcommand{\eprop}{\end{prop}}
\newcommand{\btheo}{\begin{theo}}
\newcommand{\etheo}{\end{theo}}
\newcommand{\blemm}{\begin{lemm}}
\newcommand{\brema}{\begin{rema}}
\newcommand{\erema}{\end{rema}}
\newcommand{\bexer}{\begin{exem}}
\newcommand{\eexer}{\end{exem}}
\newcommand{\bexems}{\begin{exems}}
\newcommand{\eexems}{\end{exems}}
\newcommand{\bconj}{\begin{conj}}
\newcommand{\econj}{\end{conj}}
\newcommand{\elemm}{\end{lemm}}
\newcommand{\bcoro}{\begin{coro}}
\newcommand{\ecoro}{\end{coro}}
\newcommand{\dem}{\noindent{\bf Proof. }}
\newcommand{\rem}{\noindent{\bf Remark. }}

\newcommand{\B}{{\cal B}}

\newcommand{\N}{{\cal N}}
\newcommand{\G}{{\cal G}}

\renewcommand{\O}{{\cal O}}
\newcommand{\C}{{\cal C}}

\renewcommand{\S}{{\cal S}}

\newcommand{\maths}[1]{{\mathbb #1}}  

\newcommand{\RR}{\maths{R}}
\newcommand{\NN}{\maths{N}}

\newcommand{\HH}{\maths{H}}
\newcommand{\FF}{\maths{F}}

\newcommand{\ZZ}{\maths{Z}}

\newcommand{\TT}{\maths{T}}

\newcommand{\ra}{\rightarrow}

\newcommand{\wt}[1]{{\widetilde{#1}}}
\newcommand{\wh}[1]{{\widehat{#1}}}

\newcommand{\ga}{\gamma}
\newcommand{\Ga}{\Gamma}

\newcommand{\cqfd}{\hfill$\Box$}

\newcounter{fig}

\title{On the almost sure spiraling of geodesics 
\\  in negatively curved manifolds}
 
\author{Sa'ar Hersonsky \and Fr\'ed\'eric Paulin}

\date{\today} 
 
\begin{document} 
\maketitle

\begin{abstract} 
\noindent   
Given a negatively curved geodesic metric space $M$, we study the
statistical asymptotic penetration behavior of (locally) geodesic
lines of $M$ in small neighborhoods of points, of closed geodesics,
and of other compact (locally) convex subsets of $M$. We prove
Khintchine-type and logarithm law-type results for the spiraling of
geodesic lines around these objets. As a consequence in the tree
setting, we obtain Diophantine approximation results of elements of
non-archimedian local fields by quadratic irrational ones.
\footnote{{\bf Keywords:} geodesic flow, negative curvature,
  spiraling, Khintchine theorem, logarithm law, Diophantine
  approximation, quadratic irrational.~~{\bf AMS codes:} 53 C
  22, 37 D 40, 11 J 61, 30 F 40, 37 A 45, 11 J 83}
\end{abstract}

\section{Introduction} 
\label{sec:intro}

Let $M$ be a compact connected Riemannian manifold with negative
sectional curvature. Endow the total space of the unit tangent bundle
$\pi:T^1M\ra M$ with the Bowen-Margulis measure $\mu$, which is the
maximal entropy probability measure for the geodesic flow
$(\phi_t)_{t\in\RR}$ on $T^1M$. Let $h$ be the topological entropy of
$(\phi_t)_{t\in\RR}$. In this paper, we study the
statistical asymptotic penetration behavior of (locally) geodesic
lines in various objets in $M$, as tubular neighborhoods of closed
geodesic, tubular neighborhoods of compact embedded totally geodesic
submanifolds, and other convex subsets. In this introduction, we fix a
Lipschitz map $g:\RR_+\;\ra\;\RR_+$.

We first consider a closed geodesic $C$ in $M$, and study the
spiraling of geodesics around $C$. As the geodesic flow is ergodic
with respect to $\mu$, almost every orbit in $T^1M$ is dense. Two 
geodesic lines, having at some time their unit tangent vectors
close, follow themselves closely a long time. Hence almost every 
geodesic line will
stay for arbitrarily long periods of times in a given small tubular
neighborhood of $C$. In this paper, we make this behaviour quantitative.
For that, we prove a Khintchine-type theorem, and a logarithm law-type
corollary, for geodesic lines spiraling around $C$.  Fix a small
enough $\epsilon>0$, and let $\N_\epsilon C$ be the (closed)
$\epsilon$-neighborhood of $C$.

\btheo \label{theo:introspiralclosgeod} %
If $\int_{1}^{+\infty} e^{-h\, g(t)}\;dt$ converges (resp.~diverges),
then for $\mu$-almost no (resp.~every) $v\in T^1M$, there exist
positive times $(t_n)_{n\in \NN}$ converging to $+\infty$ such that
$\pi\circ\phi_t(v)$ belongs to $\N_\epsilon C$ for every $t$ in
$[t_n,t_n+ g(t_n)]$.  \etheo


Define the penetration map ${\mathfrak p}:T^1M\times\RR \ra
[0,+\infty]$ in $\N_\epsilon C$ by ${\mathfrak
  p}(v,t)=0$ if $\pi\circ\phi_t(v)\notin \N_\epsilon C$, and otherwise
${\mathfrak p}(v,t)$ is the maximal length of an interval $I$ in $\RR$
containing $t$ such that $\pi\circ\phi_s(v)\in \N_\epsilon C$ for
every $s$ in $I$. We refer to \cite{PP} for (many) other ways to
measure the penetration of a geodesic line in the
$\epsilon$-neighborhood of $C$.

\bcoro\label{theo:introloglaw}  %
For $\mu$-almost every $v\in T^1M$,
$$
\limsup_{t\ra+\infty} \;\frac{{\mathfrak p}(v,t)}{\log t}
=\frac{1}{h}\;.
$$
\ecoro

When $M$ has constant curvature, and after a geometric translation,
Theorem \ref{theo:introspiralclosgeod} and Corollary
\ref{theo:introloglaw} follow from known results (see for instance
\cite{DMPV}, as well as the recent \cite{BV}, where the methods 
are very different).

We also prove a Khintchine-type theorem for geodesic lines spiraling
around totally geodesic submanifolds.  For the sake of simplicity in
this introduction, we only formulate it for real hyperbolic manifolds, see
Theorem \ref{theo:mainmain} for a more general statement.

\btheo\label{theo:introspiraltotgeod} %
Assume furthermore that $M$ is a real hyperbolic $n$-manifold, and 
$C$ a closed
embedded totally geodesic submanifold of dimension $k\geq 1$. Let
$\epsilon>0$ be small enough.

If $\int_{1}^{+\infty} e^{-(n-k) g(t)}\;dt$ converges
(resp.~diverges), then for $\mu$-almost no (resp.~every) $v\in T^1M$,
there exist positive times $(t_n)_{n\in \NN}$ converging to $+\infty$
such that $\pi\circ\phi_t(v)$ belongs to $\N_\epsilon C$ for every $t$
in $[t_n,t_n+ g(t_n)]$.  
\etheo

Besides totally geodesic submanifolds, one could also measure the
asymptotic spiraling of geodesic lines around other convex subsets, in
particular in hyperbolic $3$-manifolds. Recall that a subgroup $\Ga_0$
of a Kleinian group $\Ga$, whose limit set is denoted by
$\Lambda\Ga_0$, is {\it precisely invariant} if the intersection
$\ga\Lambda\Ga_0\cap\Lambda\Ga_0$ is empty for every
$\ga\in\Ga-\Ga_0$.

\btheo\label{theo:introkleinprecinv} %
Assume furthermore that $M=\Ga\backslash\HH^3_\RR$ is 
an hyperbolic $3$-manifold,
and let $\Ga_0$ be a precisely invariant quasi-fuschian subgroup of
$\Ga$. Let $\delta_0$ be the Hausdorff dimension of the limit set
$\Lambda\Ga_0$, and let $C_0$ be the image in $M$ of the convex hull of
$\Lambda\Ga_0$ in $\HH^3_\RR$. Let $\epsilon>0$ be small enough.

If $\int_{1}^{+\infty} e^{-(2-\delta_0) g(t)}\;dt$ converges
(resp.~diverges), then for $\mu$-almost no (resp.~every) $v\in T^1M$,
there exist positive times $(t_n)_{n\in \NN}$ converging to $+\infty$
such that $\pi\circ\phi_t(v)$ belongs to $\N_\epsilon C_0$ for every
$t$ in $[t_n,t_n+ g(t_n)]$.  
\etheo

All these results follow from our main result, Theorem
\ref{theo:mainmain}, which is much stronger than the above ones. We do
not require $M$ to be compact. The first two statements above remain
valid when $M$ is complete, non elementary, with a negative upper
bound on its sectional curvature, up to replacing $h$ by the critical
exponent $\delta$ of the fundamental group $\Gamma$ of $M$ acting on a
universal Riemannian covering $\wt M$ of $M$, and under some
assumptions on $\Ga$. Under these hypotheses on $M$, Theorem
\ref{theo:introkleinprecinv} is still valid, up to replacing
$\HH^3_\RR$ by $\wt M$, $\Ga_0$ by a malnormal infinite index
convex-cocompact subgroup of $\Gamma$ with crititical exponent
$\delta_0$, and $2-\delta_0$ by $\delta-\delta_0$.

Furthermore, we do not need $\wt M$ to be a manifold, Theorem
\ref{theo:mainmain} is valid in general ${\rm CAT}(-1)$ spaces, for
instance in hyperbolic buildings (i.e.~Tits buildings modeled on an
hyperbolic reflection group, see \cite{Bou2,GP,HaP} for examples).
Corollary \ref{coro:casimmhyp} is an example of an application.  In
this introduction, we quote two results in the setting of trees. The
first one will be proved in Section \ref{sec:spiralgeod}.  Let
$E[\cdot]$ be the integer part map.

\bprop\label{prop:introtree} %
Let $T$ be a locally finite tree, and ${\rm Aut}(T)$ be its locally
compact group of automorphisms, such that ${\rm Aut}(T)\backslash T$
is finite. Let $\Ga$ be a lattice in ${\rm Aut}(T)$ acting without
inversion on $T$. Let $\Ga\backslash\G T$ be the quotient by $\Ga$ of
the space $\G T$ of isometric maps $\ell:\RR\ra T$, endowed with its
geodesic flow $(\phi_t)_{t\in\RR}$ (the action of $\RR$ by
translations at the source). Let $\mu$ be the maximal entropy measure
for the action of $(\phi_t)_{t\in\RR}$ on $\Ga\backslash\G T$, and $h$
be its topological entropy.  Let $C$ be a cycle in the graph
$\Ga\backslash T$ with $L$ edges.

If $\int_{1}^{+\infty} e^{-\frac{h}{L} g(t)}\;dt$ converges
(resp.~diverges), then for $\mu$-almost no (resp.~every)
$\ell\in\Ga\backslash \G T$, there exist positive times $(t_n)_{n\in
  \NN}$ converging to $+\infty$ such that the path $t\mapsto\ell(t)$,
starting at time $t_n$, turns around $C$ at least $E[g(t_n)]$ times.
\eprop

The next result (see Section \ref{sec:dioapp}) uses the Bruhat-Tits
tree of the algebraic group ${\rm SL}_2$ over the local field $\wh
K=\FF_q((X))$ of formal Laurent series in the variable $X^{-1}$ over
the finite field $\FF_q$. Let $\mu$ be a Haar measure of $\wh K$. Let
$|\cdot|_\infty$ be the absolute value of $\wh K$. Recall that an
element of $\wh K$ is {\it irrational} if it does not belong to the
subfield $\FF_q(X)$ of rational fractions over $\FF_q$, and is {\it
  quadratic} if it is a solution of a quadratic equation with
coefficients in $\FF_q(X)$. The group ${\rm SL}_2(\FF_q[X])$ acts by
homographies on the set of quadratic irrational elements of $\wh K$,
and two of these are {\it congruent} if they are in the same orbit.
For every irrational quadratic element $\alpha$ in $\wh K$, let
$h(\alpha)= |\alpha-\alpha^*|_\infty^{-1}$, where $\alpha^*$ is the
Galois conjugate of $\alpha$, be its height (see Section
\ref{sec:dioapp}). In \cite{HP5}, we proved a $0$-$1$ measure result
for the Diophantine approximation of elements of $\wh K$ by rational
ones. The following result (see Section \ref{sec:dioapp}) is an
analogous one for the Diophantine approximation of elements of $\wh K$
by quadratic irrational ones.

\btheo\label{theo:intrononarchapproxdioph} %
Let $\varphi:\RR_+\ra\RR_+^*$ be a map with $u\mapsto \log\varphi
(e^u)$ Lipschitz. If the integral $\int_{1}^{+\infty} \varphi(t)/t
\;dt$ diverges (resp.~converges), then for $\mu$-almost every $x\in
\wh K$,
$$
\liminf\;
\frac{h(\beta)}{\varphi(h(\beta))}\;|x-\beta|_\infty=0\;
({\rm resp.} =+\infty)\;,
$$ 
where the lower limit is taken over the quadratic irrational elements
$\beta$ in $\wh K$, in any (resp.~some) congruence class, with
$h(\beta)\ra +\infty$.  
\etheo

For other number theoretic applications of the results of this paper
and of \cite{PP}, we refer to \cite{HPP}.

We first start in Section \ref{sec:backnota} by recalling properties
of the ${\rm CAT}(-1)$-spaces $X$ and their spaces at infinity
$\partial_\infty X$. We introduce, for every non empty closed convex
subset $C$ of $X$, a nice new distance-like map $d_C$ on
$\partial_\infty X-\partial_\infty C$. It generalizes Gromov's
distance when $C$ is reduced to a point (see for instance \cite{Bou}),
or Hamenst\"adt's distance when $C$ is a horoball (see
\cite{Ham}\cite[Appendix]{HP1}).

In Section \ref{sec:borelcantelli}, we present the main technical tool
of this paper, a geometric avatar of the Borel-Cantelli lemma.  This
tool will also be used in Section \ref{sec:approxpoint} (which can be
read independently of Section \ref{sec:spiralgeod}) to prove other
approximation results, both of Khintchine-type and logarithm law-type,
expressing how exactly close to a given point passes almost every
geodesic line. The following result is joint work with C.~S.~Aravinda.
See Theorem \ref{theo:introapproxpoint} for a more general statement
about the approximation of given points by almost every geodesic
lines.

\btheo\label{theo:introapproxpointconstcurv} 
Assume furthermore that $M$ has constant curvature $-1$ 
and dimension $n$, and let $x_0\in M$.  Then
for Liouville-almost every $v$ in $T^1M$,
$$ 
\liminf_{t\ra +\infty} \;\;\frac{d(\pi\circ\phi_t(v),x_0)}{\log t}\;=
\frac{1}{n-1}\;.
$$  \etheo

In Section \ref{sec:spiralgeod}, we start by explaining the general
situation covering all results
\ref{theo:introspiralclosgeod}-\ref{theo:intrononarchapproxdioph}. We
prove some new estimates on the relative geometry of the convex hulls
of subgroups of a discrete group of isometries of a ${\rm CAT}(-1)$
metric space.  Among the new geometric informations (possibly useful
for later applications), we give in Theorem
\ref{theo:estimatvolepsilneigh} a fluctuating density result. It
explains the variation in $\epsilon$ of the mass for a conformal
density of the $\epsilon$-neighborhood of the limit set of a
convex-cocompact subgroup. We then prove our main result, Theorem
\ref{theo:mainmain}.

Khintchine-type theorems and their logarithm law-type corollaries for
the spiraling of geodesic lines around cusps were obtained by
D.~Sullivan \cite{Sul} when $M$ has finite volume and constant
curvature, by D.~Kleinbock and G.~Margulis \cite{KM} if $M$ is a
finite volume locally symmetric space, by B.~Stratmann and
S.~L.~Velani \cite{SV} (see also for instance \cite{DMPV,BV}) if $M$
is geometrically finite with constant curvature, and by the authors
\cite{HP2} if $M$ is geometrically finite with variable curvature. In
this paper, which requires many new geometric inputs, our intellectual
debt to D.~Sullivan's work is still important.

\medskip
\noindent {\small {\it Acknowledgements.} The second author
  acknowledges the support of the University of Georgia at Athens for
  fruitful visits.  Theorem \ref{theo:introapproxpoint}, and most of
  the content of Section \ref{sec:approxpoint}, was essentially proved
  by Aravinda and the second author during a stay of Aravinda at the
  Universit\'e d'Orsay in 2000. We also thank V.~Kleptsyn and
  M.~Pollicot for their comments, and F.~Haglund for Remark 
  \ref{rem:haglund}.
}

\section{On convexity properties of CAT(-1)-spaces and
their discrete subgroups}
\label{sec:backnota}

We refer for instance to \cite{Bou,BH} for the definitions and basic
properties of CAT$(-1)$ metric spaces, their horospheres and their
discrete groups of isometries.  The new result in this section is the
contruction of the distance-like map $d_C$ for a convex subset $C$.

\medskip %
Let $X$ be a proper CAT$(-1)$ geodesic metric space. Its boundary at
infinity is denoted by $\partial_\infty X$. The space of geodesic
lines $\ell:\RR\ra X$ in $X$, with the compact-open topology, is
denoted by $\G X$. The geodesic flow on $\G X$ is the action of $\RR$
by translation at the source.  For every point $x_0$ in $X$, the space
$\G X$ identifies with $((\partial_\infty X\times \partial_\infty X) -
\Delta)\times \RR$, where $\Delta$ is the diagonal in $\partial_\infty
X\times \partial_\infty X$, by the map which associates to a geodesic
lines $\ell$ the triple $(\ell(-\infty), \ell(+\infty),t)$ of the
points at infinity of $\ell$ and the algebraic distance $t$ on $\ell$
(oriented from $\ell(-\infty)$ to $\ell(+\infty)$) between $\ell(0)$
and the closest point of $\ell$ to $x_0$. Given another base point
$x'_0$, this parametrization (called Hopf's) differs from the one
defined by $x'_0$ only by an additive constant on the third factor.
When $X$ is a Riemannian manifold, the map $T^1X\ra\G X$ which
associates to $v\in T^1X$ the geodesic line $t\mapsto \pi\circ\phi_t
(v)$, with $\pi:T^1X\ra X$ the canonical projection, is an
homeomorphism equivariant with respect to the actions of the isometry
group of $X$ and of the geodesic flows on $T^1X$ and $\G X$.

For every $\xi$ in $\partial_\infty X$, the {\it Buseman function at
  $\xi$} is the map $\beta_\xi$ from $X\times X$ to $\RR$ defined by
$$
\beta_\xi(x,y)=\lim_{t\ra+\infty} d(x,\xi_t)-d(y,\xi_t)\;,
$$
for any geodesic ray $t\mapsto \xi_t$ ending at $\xi$. For every
subset $A$ of $X$ and every point $x$ in $X\cup \partial_\infty X$, we
denote by $\O_{x}A$ the {\it shadow of $A$ seen from $x$}, i.e.~the
set of points at infinity of the geodesic rays or lines starting from
$x$ and meeting $A$.

\medskip
The following result, needed only in Section \ref{sec:approxpoint}, is
Proposition 3.1 of \cite{HP} (whose proof of the left inclusion is
valid under the only assumptions below).

\blemm \label{lem:HPcrelle} Let $\rho$ be a geodesic ray in $X$, with
$x=\rho(0)$ and $\xi=\rho(+\infty)$.

(1) For every $c,t>0$, the ball $B_{d_x}(\xi,ce^{-t})$ for the visual
distance (see Example 1 below) $d_x$ on $\partial_\infty X$ is
contained in the shadow $\O_x\big(B(\rho(t),c)\big)$.

(2) If $X$ is a Riemannian manifold with sectional curvature $-a^2\leq
K\leq -1$, where $a\geq 1$, then there exists
$\kappa_1,\kappa_2,\kappa_3>0$ such that for every
$c\in\;]\,0,\kappa_3]$ and every $t\geq \kappa_2$, the shadow
    $\O_x\big(B(\rho(t),c)\big)$ is contained in the ball
    $B_{d_x}(\xi,\kappa_1 \;c^{\frac{1}{a}}\;e^{-t})$. \cqfd 
\elemm

Given a point at infinity $\xi\in \partial_\infty X$ and a horosphere
$H$ centered at $\xi$, let $d_{\xi,H},d'_{\xi,H}:(\partial_\infty
X-\{\xi\})^2\ra \RR$ be the following maps. Let $\eta,\eta' \in
\partial_\infty X-\{\xi\}$.  Let $t\mapsto \eta_t$ and $t\mapsto
\eta'_t$ be the geodesic lines starting from $\xi$, passing at time
$t=0$ in $H$, and converging to $\eta$ and $\eta'$, respectively.
Define the {\it Hamenst\"adt distance} (see
\cite{Ham}\cite[Appendix]{HP1})
$$
d_{\xi,H}(\eta,\eta')=\lim_{t\ra+\infty}
e^{\frac{1}{2}d(\eta_t,\eta'_t)-t}\;,
$$
which is a distance inducing the original topology on $\partial_\infty
X-\{\xi\}$. The {\it cuspidal distance} (see \cite{HP3}) $d'_{\xi,H}$
is defined as follows: If $\eta\neq \eta'$, then 
$-\log \big(2\;d'_{\xi,H}(\eta,\eta')\big)$ is
the signed distance, along the geodesic line $]\xi,\eta[$ oriented 
from $\xi$ to $\eta$, from $H$ to the horosphere centered at $\eta$ and
meeting the geodesic line $]\xi, \eta'[$ in exactly one point. Though
not necessarily an actual distance, $d'_{\xi,H}$ is equivalent to the
Hamenst\"adt distance (see \cite[Rem.~2.6]{HP3}).

\medskip %
Let $C$ be non empty closed convex subset of $X$. (Recall that a
subset $C$ in a CAT$(-1)$ metric space is {\it convex} if $C$ contains
the geodesic segment between any two points in $C$.)  We denote by
$\partial_\infty C$ its set of points at infinity, and by $\partial C$
its boundary in $X$.  For every $\xi$ in $X\cup\partial_\infty X$, we
define {\it the closest point to $\xi$ on the convex set $C$}, denoted
by $\pi_C(\xi)$, to be the following point $p$ in $C\cup
\partial_\infty C$: If $\xi\in X$, then $p$ belongs to $C$ and
minimizes the distance between $x$ and a point of $C$; if
$\xi\in\partial_\infty X-\partial_\infty C$, then the (closed)
horoball centered at $\xi$ whose horosphere contains $p$ meets $C$
exactly at $p$; if $\xi\in\partial_\infty C$, then we define $p=\xi$.
This $p$ exists, is unique, and $\pi_C:X\cup\partial_\infty X\ra C\cup
\partial_\infty C$ is continuous, by the properties of
CAT$(-1)$-spaces. When $X$ is a Riemannian manifold, $\pi_C$ is
(outside $\partial_\infty C$) the orthogonal projection on $C$.

Let us define a distance-like map $d_C$ on $\partial_\infty
X-\partial_\infty C$. For every $\xi,\eta$ in $\partial_\infty
X-\partial_\infty C$, let $t\mapsto \xi_t$ and $t\mapsto \eta_t$ be
geodesic rays, starting at time $0$ from the closest points to $\xi$
and $\eta$ on $C$, and converging to $\xi$ and $\eta$ respectively.
Define
$$
d_C(\xi,\eta)=\lim_{t\ra+\infty} e^{\frac{1}{2}d(\xi_t,\eta_t)-t}=
\lim_{t\ra+\infty} e^{\frac{1}{2}\big(d(\xi_t,\eta_t)-d(\xi_t,\pi_C(\xi))-
d(\eta_t,\pi_C(\eta))\big)}\;.
$$
Note that these limits exist, and the second equality holds for every
geodesic rays $t\mapsto \xi_t$ and $t\mapsto \eta_t$ converging to
$\xi$ and $\eta$, respectively. For every isometry $\ga$ of $X$, we
have
$$
d_{\ga C}(\ga \xi,\ga \eta)=d_C(\xi,\eta)\;.
$$ 
In particular, any isometry of $X$ preserving $C$ preserves $d_C$. For
every $\epsilon>0$, the (closed) $\epsilon$-neighborhood $\N_\epsilon
C$ of $C$ in $X$ is still convex with $\partial_\infty(\N_\epsilon
C)=\partial_\infty C$, and
$$
d_{\N_\epsilon C}(\xi,\eta)=e^{\epsilon}\;d_C(\xi,\eta)\;.
$$

\medskip
\noindent {\bf Examples. } 

(1) If $C$ is reduced to a point $x$ in $X$, then $d_C=d_x$ is the
usual visual distance seen from $x$ on $\partial_\infty X$ (see for
instance \cite{Bou}).

(2) If $C$ is a (closed) ball of center $x$ and radius $r$, then
$d_C=e^{r}\,d_x$, as $C=\N_r \{x\}$.

(3) If $C$ is a (closed) horoball with point at infinity $\xi$ and
boundary horosphere $H$, then $d_C=d_{\xi,\, H}$ is the Hamenst\"adt
distance on $\partial_\infty X-\{\xi\}$ as recalled above. If
$(C_n)_{n\in\NN}$ is a sequence of balls converging uniformly on
compact subsets of $X$ to an horoball $C$, then the maps $d_{C_n}$
converge uniformly on compact subsets of $\partial X-\{\xi\}$ to
$d_C$.

(4) If $X$ is a metric tree, then it is easy to prove that, for every
$\xi,\eta$ in $\partial_\infty X-\partial_\infty C$,
\begin{equation}\label{eq:calcdcarb}
d_C(\xi,\eta)=\left\{\begin{array}{l}\! 
e^{\frac{1}{2}\,d(\pi_C(\xi),\,\pi_C(\eta))} 
\;\;{\rm if} \;\pi_C(\xi)\neq \pi_C(\eta)\\
\!e^{-d(\pi_C(\xi),\,q)} \;\;{\rm if} \;\pi_C(\xi)= \pi_C(\eta)\;
{\rm and} \;[\pi_C(\xi),\xi[\;\cap \,[\pi_C(\eta),\eta[ \;= 
[\pi_C(\xi),q]\;.
\end{array}\right.
\end{equation}
In particular, in a small enough neighborhood of any point $\xi_0$ in
$\partial_\infty X -\partial_\infty C$, the map $d_C$ then coincides
with the visual distance $d_{\pi_C(\xi_0)}$.

(5) Let $X$ be the real hyperbolic $n$-space $\HH^n_\RR$, and let $C$
be a complete totally geodesic submanifold of dimension $k$ with
$0<k<n$.  Let $x_0$ be a point in $C$, and $S_{x_0}(C^\perp)$ be the
sphere of unit tangent vectors at $x_0$ that are perpendicular to $C$,
endowed with the angular distance $(v,v')\mapsto\angle_{x_0}(v,v')$.
Note that the standard Euclidean distance on $S_{x_0}(C^\perp)$ is
given by $(v,v')\mapsto 2\sin\frac{\angle_{x_0}(v,v')}{2}$.  For every
$\xi\in\partial _\infty X-\partial_\infty C$, let $\pi'_C(\xi)$ be the
parallel transport to $x_0$, along a geodesic line through
$x_0,\pi_C(\xi)$, of the unit tangent vector at $\pi_C(\xi)$ of the
geodesic ray $[\pi_C(\xi),\xi[$. We thus get a map $\pi'_C:\partial
_\infty X-\partial_\infty C\ra S_{x_0}(C^\perp)$. In particular,
$(\pi_C,\pi'_C):\partial _\infty X-\partial_\infty C\ra C\times
S_{x_0}(C^\perp)$ is an homeomorphism.

\blemm \label{lem:calcdChypreel} %
For every $\xi,\eta$ in $\partial_\infty X-\partial_\infty C$,
$$
d_C(\xi,\eta)=\sqrt{\;\sinh{}^2
  \;\frac{d\big(\pi_C(\xi),\pi_C(\eta)\big)}{2}+
  \sin{}^2\;\frac{\angle_{x_0}\,
\big(\pi'_C(\xi),\pi'_C(\eta)\big)}{2}}\;.
$$
\elemm

\dem Let $\rho=d(\pi_C(\xi),\pi_C(\eta))$ and let
$\theta=\angle_{x_0}(\pi'_C(\xi),\pi'_C(\eta))$. We have to prove that
$d_C(\xi,\eta)=\frac12\sqrt{e^\rho+e^{-\rho}-2\cos\theta}$.  This last
formula follows from an easy computation using the picture below.
Recall that $\sinh b=1/\tan\alpha$, where $b$ is the hyperbolic length
of the arc of any half-circle perpendicular to the horizontal plane
between the angles $\alpha$ and $\frac{\pi}{2}$ in the upper halfspace
model of $\HH^3_\RR$ (see \cite[page 145]{Bea}).

\medskip
\noindent
\begin{minipage}{7.6cm}
  Take a copy of $\HH^3_\RR$ containing $\xi,\eta$ and a geodesic line
  passing through $\pi_C(\xi),\pi_C(\eta)$. Use the upper halfspace
  model of $\HH^3_\RR$ where this geodesic line is a vertical line
  between $0$ and $\infty$, with $\pi_C(\eta)$ above $\pi_C(\xi)$.
  Scale such that the Euclidean distance between $0$ and $\xi$ is $1$.
  Consider the points $\eta'_s,\eta_s,\xi'_s,\xi_s$ at Euclidean
  height $s$ close to $0$ on respectively $[\pi_C(\eta),\eta[\,$,
  $]\eta,\xi[$ close to $\eta$, $[\pi_C(\xi),\xi[\,$, $]\eta,\xi[$
  close to $\xi$, so that $d_C(\xi,\eta)$ is equal to
$$
\lim_{s\ra 0} \;e^{\frac{1}{2} \big(d(\xi_s,\eta_s)-
  d(\xi'_s,\pi_C(\xi))- d(\eta'_s,\pi_C(\eta))\big)}\;.
$$ 
Now just use several times the previously mentionned formula $\sinh
b=1/\tan\alpha$.  \cqfd
\end{minipage}
\begin{minipage}{6.3cm}
\begin{center}
\begin{picture}(0,0)%
\includegraphics{fig_calcdistlike.pstex}%
\end{picture}%
\setlength{\unitlength}{3729sp}%
\begingroup\makeatletter\ifx\SetFigFont\undefined%
\gdef\SetFigFont#1#2#3#4#5{%
  \reset@font\fontsize{#1}{#2pt}%
  \fontfamily{#3}\fontseries{#4}\fontshape{#5}%
  \selectfont}%
\fi\endgroup%
\begin{picture}(2877,2590)(811,-3980)
\put(2251,-1546){\makebox(0,0)[lb]{\smash{{\SetFigFont{11}{13.2}{\rmdefault}{\mddefault}{\updefault}{\color[rgb]{0,0,0}$C\cap\HH^3_\RR$}%
}}}}
\put(2813,-2191){\makebox(0,0)[lb]{\smash{{\SetFigFont{11}{13.2}{\rmdefault}{\mddefault}{\updefault}{\color[rgb]{0,0,0}$\rho$}%
}}}}
\put(2476,-3376){\makebox(0,0)[lb]{\smash{{\SetFigFont{11}{13.2}{\rmdefault}{\mddefault}{\updefault}{\color[rgb]{0,0,0}$\theta$}%
}}}}
\put(2686,-3676){\makebox(0,0)[lb]{\smash{{\SetFigFont{11}{13.2}{\rmdefault}{\mddefault}{\updefault}{\color[rgb]{0,0,0}$\sqrt{e^{2\rho}+1-2e^\rho\cos\theta}$}%
}}}}
\put(3354,-3301){\makebox(0,0)[lb]{\smash{{\SetFigFont{11}{13.2}{\rmdefault}{\mddefault}{\updefault}{\color[rgb]{0,0,0}$\xi$}%
}}}}
\put(2604,-1905){\makebox(0,0)[lb]{\smash{{\SetFigFont{11}{13.2}{\rmdefault}{\mddefault}{\updefault}{\color[rgb]{0,0,0}$\pi_C(\eta)$}%
}}}}
\put(811,-3841){\makebox(0,0)[lb]{\smash{{\SetFigFont{11}{13.2}{\rmdefault}{\mddefault}{\updefault}{\color[rgb]{0,0,0}$s$}%
}}}}
\put(1891,-3886){\makebox(0,0)[lb]{\smash{{\SetFigFont{11}{13.2}{\rmdefault}{\mddefault}{\updefault}{\color[rgb]{0,0,0}$\eta$}%
}}}}
\put(2791,-2941){\makebox(0,0)[lb]{\smash{{\SetFigFont{11}{13.2}{\rmdefault}{\mddefault}{\updefault}{\color[rgb]{0,0,0}$1$}%
}}}}
\put(1666,-3526){\makebox(0,0)[lb]{\smash{{\SetFigFont{11}{13.2}{\rmdefault}{\mddefault}{\updefault}{\color[rgb]{0,0,0}$\eta'_s$}%
}}}}
\put(3061,-2941){\makebox(0,0)[lb]{\smash{{\SetFigFont{11}{13.2}{\rmdefault}{\mddefault}{\updefault}{\color[rgb]{0,0,0}$\xi_s$}%
}}}}
\put(1973,-3399){\makebox(0,0)[lb]{\smash{{\SetFigFont{11}{13.2}{\rmdefault}{\mddefault}{\updefault}{\color[rgb]{0,0,0}$\eta_s$}%
}}}}
\put(3346,-2836){\makebox(0,0)[lb]{\smash{{\SetFigFont{11}{13.2}{\rmdefault}{\mddefault}{\updefault}{\color[rgb]{0,0,0}$\xi'_s$}%
}}}}
\put(2008,-2356){\makebox(0,0)[lb]{\smash{{\SetFigFont{11}{13.2}{\rmdefault}{\mddefault}{\updefault}{\color[rgb]{0,0,0}$\pi_C(\xi)$}%
}}}}
\put(2138,-3196){\makebox(0,0)[lb]{\smash{{\SetFigFont{11}{13.2}{\rmdefault}{\mddefault}{\updefault}{\color[rgb]{0,0,0}$e^\rho$}%
}}}}
\put(2416,-3038){\makebox(0,0)[lb]{\smash{{\SetFigFont{11}{13.2}{\rmdefault}{\mddefault}{\updefault}{\color[rgb]{0,0,0}$0$}%
}}}}
\end{picture}%

\end{center}
\end{minipage}

\medskip 
In particular, if $X=\HH^2_\RR$, if $C$ is a geodesic line and if
$\xi,\eta$ are in the same component of $\partial_\infty
X-\partial_\infty C$, then
$$
d_C(\xi,\eta)=\sinh \;\frac{d(\pi_C(\xi),\pi_C(\eta))}{2}\;.
$$
By taking $a,b,c$ in the same component of $\partial_\infty X-
\partial_\infty C$ such that $d(\pi_C(a),\pi_C(b))=d(\pi_C(b),\pi_C(c))=
\frac12 d(\pi_C(a),\pi_C(c))$ are big enough, we see that $d_C$ does not
satisfy the triangle inequality, hence is not a distance.

\medskip
After these examples, let us go back to the general situation on $X,C$, 
and let us prove some results saying that at least on compact subsets, 
the map $d_C$ behaves quite like a distance.

\blemm \label{lem:propdistconv} (1) For every $x_0$ in $X$, for every
compact subset $K$ of $\partial_\infty X-\partial_\infty C$, there
exists a constant $c_K>0$ such that for every $\xi,\eta$ in $K$, we
have
$$
\frac{1}{c_K}d_{x_0}(\xi,\eta)\leq d_C(\xi,\eta)\leq 
c_K\,d_{x_0}(\xi,\eta)\;.
$$
(2) For every $\xi$ in $\partial_\infty X-\partial_\infty C$, the map
$\eta\mapsto d_C(\xi,\eta)$ is proper on $\partial_\infty
X-\partial_\infty C$.

\noindent (3) For every $\xi,\eta$ in $\partial_\infty X -
\partial_\infty C$,
$$
(3-2\sqrt{2})\;e^{\frac{1}{2}d(\pi_C(\xi),\pi_C(\eta))}\;
e^{-d(C,\;]\xi,\eta[\,)}\leq d_C(\xi,\eta)\leq 
e^{\frac{1}{2}d(\pi_C(\xi),\pi_C(\eta))}\;.
$$

\noindent (4) There exist universal constants $c,c'>0$ such that for
every $\xi,\eta$ in $\partial_\infty X-\partial_\infty C$, if
$d_{C}(\xi,\eta)\leq c$, then $C$ and the geodesic line $]\xi,\eta[$
are disjoint, and
$$
\frac{1}{c'}\;e^{-d(C,\;]\xi,\eta[\,)}\;\leq d_C(\xi,\eta)\leq 
\;c'\;e^{-d(C,\;]\xi,\eta[\,)}\;.
$$
\elemm

Note that by hyperbolicity, $\min \{d(\pi_C(\xi),\pi_C(\eta)),
d(C,]\xi, \eta[\,)\}$ is, for every $\xi,\eta$, less than a universal
constant.

\medskip \dem 
For every $x_0$ in the convex subset $C$ and $\xi$ in $\partial_\infty X-
\partial_\infty C$, by the triangle inequality and the ${\rm CAT}(-1)$
inequality, we have
$$
d(\xi_t,x_0)\leq d(\xi_t,\pi_C(\xi))+d(x_0,\pi_C(\xi))\leq
d(\xi_t,x_0) +2\log(1+\sqrt{2})\;,
$$ 
with $\xi_t$ as above. Hence for every $\xi,\eta$ in $\partial_\infty
 X-\partial_\infty C$,
 \begin{equation}\label{eq:reldcdx0}
(3-2\sqrt{2}) \;d_{x_0}(\xi,\eta)\leq d_C(\xi,\eta)\;
e^{-\frac{1}{2}\big(d(x_0,\pi_C(\xi))+d(x_0,\pi_C(\eta))\big)}\leq
d_{x_0}(\xi,\eta)\;.
\end{equation}
The first result easily follows. By taking $x_0=\pi_C(\xi)$ in the
lower bound of Equation \eqref{eq:reldcdx0}, the second assertion also
follows. By taking $x_0$ to be the midpoint of the geodesic segment
$[\pi_C(\xi),\pi_C(\eta)]$ in Equation \eqref{eq:reldcdx0}, and since
$d_{x_0}\leq 1$, the upper bound in the third assertion follows.

By the triangle inequality, $d_{x_0}(\xi,\eta)\geq
e^{-d(x_0,\;]\xi,\eta[\,)}$. Hence, by taking $x_0$ in Equation
\eqref{eq:reldcdx0} to be the closest point of $C$ to $]\xi,\eta[$ if
$C$ and $]\xi,\eta[$ are disjoint, or any point in $C\cap\,]\xi,\eta[$
otherwise, and by using again the triangle inequality, the lower bound
in the third assertion follows.

The last assertion follows by standard techniques of approximation by
trees (see for example \cite[page 33]{GH}). \cqfd

\medskip In particular, the non negative symmetric map $d_C$ vanishes
on and only on the diagonal of $(\partial_\infty X-\partial_\infty
C)^2$.  But as seen above, $d_C$ is not always a distance.

It also follows from Lemma \ref{lem:propdistconv} (1) that the uniform
structure (see for instance \cite{Bourb}) defined (on compacts
subsets) by the family $\left(\{(x,y)\in (\partial_\infty
  X-\partial_\infty C)^2 \;:\;d_C(\xi,\eta)\leq \epsilon\}
\right)_{\epsilon>0}$ is isomorphic (on compacts subsets) to the
uniform structure defined by the distance $d_{x_0}$.

\medskip \rem Though we won't need it in this paper, here is a formula
expressing the distance-like map $d_C$, when $C=L$ is a geodesic line
with endpoints $L_-, L_+$, in terms of the Hamenst\"adt distance and
the cuspidal distance: For every $\xi,\eta$ in $\partial_\infty
X-\{L_-,L_+\}$, for every horosphere $H$ small enough centered at
$L_-$,
$$
d_C(\xi,\eta)=\frac{d_{L_-,H}(\xi,\eta)}{2\big(d'_{L_-,H}(\xi,L_+)
d'_{L_-,H}(\eta,L_+)\big)^{\frac{1}{2}}}\;.
$$

\dem Let $H_\xi$ (resp.~$H_\eta$) be the horosphere centered at $\xi$
(resp.~$\eta$) passing through $\pi_C(\xi)$ (resp.~$\pi_C(\eta)$). Let
$h_\xi$ (resp.~$h_\eta$) be the intersection point of $H_\xi$
(resp.~$H_\eta$) with the geodesic line $]\xi,L_-[$
(resp.~$]\eta,L_-[$). Then
$$d_C(\xi,\eta)=
\lim_{t\ra+\infty}\;e^{\frac{1}{2}\big(d(\xi_t,\eta_t)-d(\xi_t,h_\xi)-
d(\eta_t,h_\eta)\big)}=d_{L_-,H}(\xi,\eta)\;
e^{\frac{1}{2}\big(d(h_\xi,H)+d(h_\eta,H)\big)}
\;,
$$
which proves the result.  
\cqfd

\bigskip 
Let $\Ga$ be a discrete group of isometries of $X$.  Its
limit set is denoted by $\Lambda\Ga$, and if $\Lambda\Ga$ contains at
least two points, then the convex hull of $\Lambda\Ga$ is denoted by
$\C \Ga$.  Recall that $\partial_\infty\C\Gamma=\Lambda\Gamma$.  The
{\it critical exponent} of $\Ga$ is the unique number $\delta_\Ga$ in
$[0,+\infty]$ such that the {\it Poincar\'e series}
$P_{x_0,\Ga}(s)=\sum_{\ga\in\Ga} e^{-s \;d(x_0,\ga x_0)}$ of $\Ga$
converges for $s>\delta_\Ga$ and diverges for $s<\delta_\Ga$, where
$x_0$ is any point in $X$. The group $\Ga$ is called {\it of divergent
  type} if its Poincar\'e series diverges at $s=\delta_\Ga$.  The
group $\Ga$ is {\it non elementary} if $\Lambda\Ga$ contains at least
three points, and we have then $\delta_\Ga>0$. Note that when $X$ is a
Riemannian manifold and $\Ga$ is torsion free with compact quotient
$X/\Ga$, then the critical exponent $\delta_\Ga$ of $\Ga$ is the
topological entropy of the geodesic flow of $X/\Ga$ (see for instance
\cite{Man}).

If $\delta\in\;]0,+\infty[$, a {\it conformal} (or Patterson-Sullivan)
{\it density of dimension $\delta$} for $\Ga$ is a family
$(\mu_x)_{x\in X}$ of finite Borel measures on $\partial_\infty X$,
such that
\begin{itemize}
\item $\forall\;\ga\in\Ga\;,\;\;\ga_*\mu_x=\mu_{\ga x}$,
\item $\forall\;x,y\in X,\;\forall\;\xi\in\partial_\infty X
  \;,\;\;\frac{d\mu_x}{d\mu_y}(\xi)=e^{-\delta \beta_\xi(x,y)}$.
\end{itemize}
Using Hopf's parametrization with respect to any base point $x_0$ of
$X$, the (Patterson-Sullivan-){\it Bowen-Margulis measure} associated
to this family is the measure $\wt{\mu}_{\mbox{\tiny BM}}$ on $\G X$
given by
$$
d\,\wt{\mu}_{\mbox{\tiny BM}}=
\frac{d\mu_{x_0}(\xi)\,d\mu_{x_0}(\eta)\,dt}
{d_{x_0}(\eta,\xi)^{2\delta}} \;.
$$
This measure on $\G X$ is independant of $x_0$, invariant by the
action of $\Ga$ and by the geodesic flow (and by the time reversal
$\ell\mapsto \{t\mapsto \ell(-t)\}$), hence defines a measure
$\mu_{\mbox{\tiny BM}}$ on $\Ga\backslash \G X$ which is 
invariant by the quotient
geodesic flow (see for instance \cite{Bou,Rob}). Note that if 
$\mu_{\mbox{\tiny BM}}$  is finite, then $\delta=\delta_\Ga$ and $\Ga$ 
is of divergent type (see \cite[page 18-19]{Rob}). 

If $\Ga$ is of divergent type with a finite non zero critical exponent
$\delta$, then (see for instance \cite{Bou}) there exists a conformal
density of dimension $\delta$ for $\Ga$, which is unique up to a
positive scalar factor, and which is ergodic with respect to the
action of $\Ga$ on $\partial_\infty X$. The Bowen-Margulis measure
associated to any such conformal family (both on $\G X$ and on
$\Ga\backslash \G X$) will be called a {\it Bowen-Margulis measure} of
$\Ga$ (it is also uniquely defined up to a positive scalar factor).
When $X$ is a manifold and $\Ga$ acts freely on $X$ with compact
quotient, the Bowen-Margulis measure on the unit tangent bundle of the
compact negatively curved manifold $M=\Ga\backslash X$, normalized to
be a probability measure, is the maximal entropy probability measure
for the geodesic flow of $M$ (via the canonical identification of $\G
X$ and $T^1 X$), see for instance \cite{Kai}. When furthermore $X$ has
constant curvature, then the Bowen-Margulis measure and the Liouville
measure (when both are normalized) coincide on $M$.

The following result, which is obvious by definition of the Bowen-Margulis
measure, will be used in the sections \ref{sec:spiralgeod} and
\ref{sec:approxpoint}.

\blemm\label{lem:bowmarabscontpatsul} %
Let $\pi_+:\G X\ra\partial_\infty X$ be the map $\ell\mapsto
\ell(+\infty)$. Let $\wt{\mu}_{\mbox{\tiny BM}}$ be the Bowen-Margulis
measure on $\G X$ associated to a conformal family $(\mu_x)_{x\in X}$
for $\Ga$.  Then the preimage by $\pi_+$ of a set of measure $0$
(resp.~$>0$) for $\mu_x$ (for some (equivalently for any) $x$ in $X$) has
measure $0$ (resp.~$>0$) for $\wt{\mu}_{\mbox{\tiny BM}}$.  \cqfd
\elemm

Besides its invariance under $\Ga$ and the geodesic flow, and
its ergodicity on $\Ga\backslash\G X$, this is the only property of
the Bowen-Margulis measure $\wt{\mu}_{\mbox{\tiny BM}}$ on $\G X$ that
will be used in this paper. In particular, we may replace
$\wt{\mu}_{\mbox{\tiny BM}}$ by any other measure satisfying these
invariance properties and this lemma, as for instance the Knieper
measure (see \cite{Kni}).

\medskip
The group $\Ga$ is said to be {\it convex-cocompact} if $\Lambda\Ga$
contains at least two points, and if the action of $\Ga$ on $\C \Ga$
has compact quotient. In particular, the group generated by an
hyperbolic isometry of $X$ is convex-cocompact, with critical exponent
$0$. In fact, if $\Ga$ is convex-cocompact, then its critical
exponent is $0$ if and only if $\Ga$ has an index two subgroup
generated by an hyperbolic isometry of $X$. Note that $\partial_\infty
\C\Ga=\Lambda\Ga$, and that  if $\Ga$ is
convex-cocompact then $\Ga$ is of divergent type 
(see for instance \cite{Bou,Rob}).

For every $f,g:\NN\ra[0,+\infty[$, write $f\asymp g$ if there exists a
constant $c\geq 1$ such that $\frac{1}{c}f \leq g \leq c f$. For every
$x_0$ in $X$, if $\Ga$ is convex-cocompact and non elementary, with
critical exponent $\delta_\Ga$, then
$$
{\rm Card}(B(x_0,n)\cap\Ga x_0)\asymp e^{\delta_\Ga n}
$$
(see for instance \cite{Rob}, where others, much more general,
assumptions on $\Ga$ are given for this property to hold. This is the
case for example when the Bowen-Margulis measure $\mu_{\mbox{\tiny
    BM}}$ of $\Ga$ is finite (and the length spectrum is non
arithmetic), see \cite[page 56]{Rob}).

\blemm\label{lem:convcocinfindexpstrinf} %
Let $\Ga_0$ be a convex-cocompact subgroup with infinite index in a
discrete group of isometries $\Ga$ of $X$. Let $\delta_0$ and $\delta$
be the critical exponents of $\Ga_0$ and $\Ga$ respectively. Then
$\delta_0<\delta$.  \elemm

\dem %
This is well-known (see for instance \cite{Fur} in a special case).
\cqfd

\medskip 
Recall that the {\it virtual normalizer} ${\rm N}\Ga_0$ of a
convex-cocompact subgroup $\Ga_0$ of $\Ga$ is the stabilizer in $\Ga$
of the limit set $\Lambda\Ga_0$. It contains the normalizer of $\Ga_0$
in $\Ga$, and it contains $\Ga_0$ with finite index (see for instance
\cite{KS,Arz}).

Recall that a subgroup $H$ of a group $G$ is {\it malnormal} if, for
every $g$ in $G-H$, we have $gHg^{-1}\cap H=\{1\}$. We will say that a
subgroup $H$ of a group $G$ is {\it almost malnormal} if, for every
$g$ in $G-H$, the subgroup $gHg^{-1}\cap H$ is finite. Note that
malnormal implies almost malnormal, and that the converse is true if
the ambient group is torsion free.

The following result is folklore, we provide a proof because we
couldn't find a precise reference.

\bprop\label{prop:equivmalnormal} Let $\Ga_0$ be a convex-cocompact
subgroup of a discrete group $\Ga$ of isometries of $X$, then the
following assertions are equivalent.
\begin{enumerate}
\item[(1)] 
  $\Ga_0$ is almost malnormal in $\Ga$;
\item[(2)] 
  the limit set of $\Ga_0$ is precisely invariant, i.e.~for
  every $\ga\in\Ga-\Ga_0$, the set $\Lambda\Ga_0\cap\ga\Lambda\Ga_0$
  is empty;
\item[(3)] 
  $\C\Ga_0\cap \ga \C\Ga_0$ is compact for every $\ga\in\Ga-\Ga_0$;
\item[(4)] for every $\epsilon>0$, there exists
  $\kappa=\kappa(\epsilon)>0$ such that ${\rm diam}\big(\N_\epsilon
  \C\Ga_0\cap \ga\N_\epsilon \C\Ga_0\big)\leq \kappa$ for every
  $\ga\in\Ga-\Ga_0$.
\end{enumerate}
\eprop

The convex hull in $X$ of the limit set of a convex-cocompact
subgroup is non compact. Hence an almost malnormal convex-cocompact
subgroup of $\Ga$ is equal to its virtual normalizer, by (3).

\medskip \dem 
As $\partial_\infty \C\Ga_0=\Lambda\Ga_0$, it is clear
that (4) implies (3), which implies (2), which implies (1).

Let us prove that (1) implies (4). Let $C_0=\C\Ga_0$ and $\epsilon>0$.
Assume by absurd that for every $n$ in $\NN$, there exists $\ga_n$ in
$\Ga-\Ga_0$ and $x_n,y_n$ in $\N_{\epsilon}C_0\, \cap\,
\ga_n\N_{\epsilon} C_0$ with $d(x_n,y_n)\geq n$.  As $\Ga_0\backslash
\N_{\epsilon}C_0$ is compact and the action of $\Ga_0$ is isometric,
there exists $R>0$ such that $\Ga_0 B(x,R)$ contains
$\N_{\epsilon}C_0$ for every $x$ in $\N_{\epsilon}C_0$.  As $\ga_n$ is
an isometry, we also have that $\ga_n \Ga_0\ga_n^{-1} B(y,R)$ contains
$\ga_n \N_{\epsilon}C_0$, for every $y$ in $\ga_n \N_{\epsilon}C_0$.
Up to conjugating $\ga_n$ by an element of $\Ga_0$, we may assume that
$x_n$ stays in a compact subset $K$ of $X$, and we define $K'=\N_RK$,
which is compact.  As $\Ga$ is discrete, the number $N$ of elements
$\ga$ in $\Ga$, such that $\ga K'\cap K'$ is non empty, is finite. As
$\Ga_0$ is convex-cocompact, the upper bound of the cardinals of the
finite subgroups of $\Ga_0$ is finite. Hence, as $\Ga_0$ is almost
malnormal, there exists $N'\in\NN$ such that for every $\ga$ in
$\Ga-\Ga_0$, the cardinal of $\ga\Ga_0\ga^{-1}\cap\Ga_0$ is at most
$N'-2$. Take $n$ in $\NN$ with $n>N N' {\rm ~diam~} K'$. Subdivide the
segment between $x_n$ and $y_n$ in points $u_0=x_n,u_1,\dots,
u_{NN'}=y_n$, such that $d(u_k,u_{k+1})>{\rm ~diam~} K'$ for $0\leq
k\leq NN'-1$. As $K'$ contains $B(x_n,R)$ and $x_n,y_n$ belong to the
convex subset $\N_{\epsilon}C_0\, \cap\, \ga_n\N_{\epsilon} C_0$, for
$0\leq k\leq NN'$, there exist $\alpha_k,\beta_k$ in $\Ga_0$ such that
$u_k\in\alpha_k K'$ and $u_k\in\ga_n\beta_k \ga_n^{-1} K'$. Note that
$\alpha_k\neq \alpha_j$ if $k\neq j$, as $d(u_k,u_j)>{\rm ~diam~} K'$.
By the definition of $N$, there exists $(k_j)_{1\leq j\leq N'}$ with
$\alpha_{k_j}^{-1}\ga_n\beta_{k_j} \ga_n^{-1}=
\alpha_{k_{1}}^{-1}\ga_n \beta_{k_{1}}\ga_n^{-1}$ for $1\leq j\leq
N'$. Hence $\ga_n\beta_{k_j}\beta_{k_{1}}^{-1}
\ga_n^{-1}=\alpha_{k_j}\alpha_{k_{1}}^{-1}$, for $2\leq j\leq N'$,
which contradicts the fact that the cardinal of
$\ga_n\Ga_0\ga_n^{-1}\cap\Ga_0$ is at most $N'-2$. \cqfd

\medskip
\rem The fact that the first two assertions are equivalent follows
also from the well-known equality
$$
\Lambda\Ga_0\cap\ga\Lambda\Ga_0=\Lambda(\Ga_0\cap\ga\Ga_0\ga^{-1})
\;,
$$ see for instance
\cite[Coro.~3]{SS} for a proof in a special case.

\section{A geometric avatar of the Borel-Cantelli lemma} 
\label{sec:borelcantelli}

The main technical tool of this paper is the following result, which
is a suitable enhancement of the Borel-Cantelli Lemma.

\btheo\label{theo:resucborelcantelli} %
Let $(Z,\mu)$ be a measured space with $\mu(Z)$ finite, and
$(B_i(\epsilon))_{i\in I,\; \epsilon\,\in\;]\,0,+\infty[}$ a family of
measurable subsets in $Z$, non-decreasing in $\epsilon$ (for the
inclusion), endowed with a map $i\mapsto n_i$ from $I$ to $\NN$ such
that $I_n=\{i\in I\;:\; n_i=n\}$ is finite for every $n$.  Let
$f_1,f_2, f_3, f_4$ be maps from $\NN$ to $]\,0,+\infty[$ and $f_5$ a
map from $]\,0,+\infty[$ to itself. Let $E$ be the (measurable) set of
points in $Z$ belonging to infinitely many subsets $B_i(f_3(n_i))$ for
$i$ in $I$.

\smallskip
\noindent $[A]$ Assume that $f_3\leq f_2$ and that there exists $c\geq
1$ such that, for every $n$ in $\NN$, $i$ in $I$ and $\epsilon\in \;
]\,0,f_2(n_i)]$, one has $ {\rm Card~} I_n \leq c f_1(n)$ and
$\mu(B_i(\epsilon)) \leq c f_4(n_i)f_5(\epsilon)$.  If the series
$\sum_{n=0}^\infty f_1(n) f_4(n) f_5(f_3(n))$ converges, then
$\mu(E)=0$.

\smallskip
\noindent $[B]$ Assume that there exists $c\geq 1$ such that 
\begin{enumerate}
\item[(1)] $f_3\leq f_2$,
\item[(2)] $\frac{1}{f_5\circ f_2}\leq f_4f_1$,
\item[(3)] there exists $c',c''>1$ such that for every
  $\epsilon,\epsilon'>0$, if $\epsilon'\leq c'\epsilon$, then
  $f_5(\epsilon')\leq c''f_5(\epsilon)$,
\item[(4)] for every $n$ in $\NN$, one has $\frac{1}{c} f_1(n) \leq
  {\rm Card~} I_n \leq c f_1(n)$,
\item[(5)] for every $i$ in $I$ and $\epsilon\in\;]\,0,f_2(n_i)]$, we
  have $\frac{1}{c} f_4(n_i)f_5(\epsilon) \leq \mu(B_i(\epsilon)) \leq
  c f_4(n_i)f_5(\epsilon)$,
\item[(6)] for every $n$ in $\NN$, the subsets $B_i(f_2(n))$ for $i$
  in $I_n$ are pairwise disjoint,
\item[(7)] for every $i,j$ in $I$ with $n_i<n_j$, if the intersection
  of $B_j(f_3(n_j))$ and $B_i(f_3(n_i))$ is non empty, then
  $B_j(f_2(n_j))$ is contained in $B_i(cf_3(n_i))$.
\end{enumerate}
If the series $\sum_{n=0}^\infty f_1(n) f_4(n) f_5(f_3(n))$ 
diverges, then $\mu(E)>0$.
\etheo

Note that (except for the convergence of the series) every hypothesis
of Case [A] is part of an hypothesis (1)-(5) of Case [B]. Hence when
checking the hypotheses when we want to apply both cases of this
theorem, we will only chek the ones of Case [B].

\medskip
\dem %
For $i$ in $I$ and $n$ in $\NN$, let $B_i=B_i(f_3(n_i))$ and
$A_n=\bigcup_{i\in I_n} B_i$, so that $E=\bigcap_{n\in\NN}
\bigcup_{k\geq n} A_k$.

Under the assumptions of [A], by the subadditivity of $\mu$, we have
the inequality $\mu(A_n)\leq c^2f_1(n) f_4(n)f_5( f_3(n))$. Therefore
the end of the proof is standard: If the series $\sum_{n=0}^\infty
f_1(n)f_4(n) f_5(f_3(n))$ converges, then the sequence
$u_k=\sum_{n=k}^\infty f_1(n)f_4(n) f_5(f_3(n))$ tends to $0$,
therefore
$$
\mu(E)=\lim_{n\ra\infty} \mu\left(\;\bigcup_{k= n}^\infty
  A_k\right)\leq \lim_{n\ra\infty} c^2u_n=0\;.
$$

\medskip %
Assume now that the assumptions of [B] hold. We first claim that
$$f_1(n)f_4(n) f_5(f_3(n))\leq c^{2}\mu(A_n)\;.\;\;\;(*)$$ 
Indeed, the balls $B_i$ for $i$ in $I_n$ are pairwise disjoint by (1)
and (6), since the subsets $B_i(r)$ are non-decreasing in $r$. By the
additivity of $\mu$, by the lower bounds in (4) and (5), the
inequality (*) hence follows.

In particular, $\sum\mu(A_n)$ diverges if $\sum f_1(n) f_4(n) f_5(
f_3(n))$ diverges.

\medskip Now, let $n,m$ be in $\NN$ with $n<m$. By the properties (6)
and (7), for every $i$ in $I_n$, we have
$$
\mu\left(B_i(cf_3(n_i))\right) \geq {\rm Card}\{j\in I_m\;:\;
B_j\cap B_i\neq \emptyset\}\;\min_{j\in
I_m}\mu\left(B_j(f_2(m))\right) \;.$$
Hence by (5)
$$
{\rm Card}\{j\in I_m\;:\; B_j\cap B_i\neq \emptyset\} \leq
\frac{c f_4(n)f_5(cf_3(n))}{\frac{1}{c}f_4(m)f_5(f_2(m))} \;. \;\;(**)
$$
Therefore
$$\begin{array}{ccl} \mu(A_n\cap A_m) &\leq &
\sum_{i\in I_n}\sum_{j\in
I_m\,,\; B_j\cap B_i\neq \emptyset} \;\mu (B_j)\medskip 
\\\medskip
&\leq & c f_1(n) \times
\frac{c f_4(n)f_5(cf_3(n))}{\frac{1}{c}f_4(m)f_5(f_2(m))} \times 
c f_4(m)f_5(f_3(m))\\\medskip 
&\leq & c^4 (c'')^{\frac{\log c}{\log c'}+1}
f_1(n)f_4(n)f_5(f_3(n))f_1(m)f_4(m)f_5(f_3(m))\\
&\leq & c^8 (c'')^{\frac{\log c}{\log c'}+1}\mu(A_n)\mu(A_m)\;.
\end{array}
$$
The second inequality follows from (4), (**) and (5), the third
inequality follows from (2) and an iterated application of (3), and
the last one from (*).

\medskip The following Borel-Cantelli Lemma is well-known (see for
instance \cite{Spr}).
 
\btheo\label{theo:borel_cantelli} Let $(Z,\nu)$ be a probability
space. Let $(A_n)_{n\in\NN}$ be a sequence of measurable subsets of
$Z$ such that there exists a constant $c>0$ with $\nu(A_n\cap A_m)\leq
c\nu(A_n)\nu(A_m)$ for every distinct integers $n,m$.  Let $ A_\infty=
\bigcap_{n\in\NN}\bigcup_{k\geq n} A_k$. Then $\nu(A_\infty)>0$ if and
only if $\displaystyle \sum_{n=0}^{\infty} \nu(A_n)$ diverges. \cqfd
\etheo

The result then follows. 
\cqfd

\section{Spiraling geodesics} 
\label{sec:spiralgeod} 


Let $X$ be a proper ${\rm CAT}(-1)$ geodesic metric space. Let $\Ga$
be a non elementary discrete group of isometries of $X$, with finite
critical exponent $\delta$. Let $\Ga_0$ be an almost malnormal
convex-cocompact subgroup of infinite index in $\Ga$ with critical
exponent $\delta_0$, and let $C_0=\C\Gamma_0$. Let $\pi_{C_0}:
X\cup\partial_\infty X\ra C_0\cup\partial_\infty C_0$ be the closest
point map.  By Lemma \ref{lem:convcocinfindexpstrinf}, the number
$\delta_0$ belongs to $[0,\delta[$. Moreover, it follows from Section
\ref{sec:backnota} that $C_0$ is non compact and that $\Ga_0$ is the
stabilizer in $\Ga$ of $C_0$.

\medskip\noindent
{\bf Examples. } 

(1) Let $\ga_0$ be an hyperbolic element of $\Ga$, let $C_0$ be its
translation axis and let $\Ga_0$ be the stabilizer of $C_0$ (which is
virtually infinite cyclic, and infinite cyclic when $\Ga$ is torsion
free). Since $\Ga$ is non elementary, the subgroup $\Ga_0$ has
infinite index. Furthermore, if $\ga\in\Ga$ and $\ga\Ga_0\ga^{-1}
\cap\Ga_0$ is infinite, then $\ga$ conjugates some hyperbolic element
of $\Ga_0$ to another one. The image by an element $\ga$ in $\Ga$ of
the translation axis of an hyperbolic element $\alpha$ of $\Ga$ is the
translation axis of $\ga\alpha\ga^{-1}$.  Hence $\ga$ preserves $C_0$,
therefore belongs to $\Ga_0$. Therefore $\Ga_0$ is an almost malnormal
convex-cocompact subgroup of infinite index in $\Ga$ with critical
exponent $\delta_0=0$.

(2) Let $M$ be a complete Riemannian manifold with dimension $n\geq 2$
and sectional curvature at most $-1$, and $\pi:X\ra M$ be a
universal Riemannian covering, with covering group $\Ga$. Let $M_0$ be
a compact connected embedded totally geodesic submanifold in $M$ of
dimension $k$ with $1\leq k\leq n-1$, let $C_0$ be a connected
component of the premimage of $M_0$ in $X$, and let $\Ga_0$ be the
stabilizer of $C_0$ in $\Ga$ (with good choices of base points, $\Ga$
can be identified with the fundamental group of $M$, and $\Ga_0$ with the
image in the fundamental group of $M$ of the fundamental group of
$M_0$). Then $\Ga_0$ is an almost malnormal (for instance by
Proposition \ref{prop:equivmalnormal} (3)) convex-cocompact subgroup
of infinite index in $\Ga$. If $M$ has constant sectional curvature
$-1$, then
$\delta=n-1$ and $\delta_0=k-1$.

(3) Let $X=\HH^3_\RR$ be the real hyperbolic space of dimension $3$,
and $\Ga$ be a  Kleinian group. If $\Ga_0$ is a
precisely invariant quasi-fuschian subgroup, without parabolic
elements, of infinite index in $\Ga$, then $\Ga_0$ is an almost
malnormal (by Proposition \ref{prop:equivmalnormal} (2))
convex-cocompact subgroup of infinite index in $\Ga$.

\bigskip %
After these examples, let us proceed. Denote by $R_0$ the set of 
double cosets
$$
R_0=\Ga_0\backslash (\Ga-\Ga_0)/\Ga_0\;.
$$
For every $r=[\ga]$ in $R_0$,  define
$$
D(r)=d(C_0,\ga C_0)\in [0,+\infty[\;,
$$
which does not depend on the representative $\ga$ of $r$. The next
result says that the subset $\{D(r)\;:\;r\in R_0\}$ of $[0,+\infty[$
is discrete, with finite multiplicities.

\blemm \label{lem:depthsdiscrete} For every $c\geq 0$, the set of
elements $r$ in $R_0$ such that $D(r)\leq c$ is finite.  
\elemm

\dem %
For every $c\geq 0$, assume that there exists a sequence of pairwise
distinct elements $([\ga_i])_{i\in \NN}$ in $R_0$ such that
$D([\ga_i])\leq c$ for every $i$.  Fix $x_*$ in $C_0$, and let $D$ be
the diameter of $\Ga_0\backslash C_0$. For every $i$ in $\NN$, let
$x_i$ in $C_0$ and $y_i$ in $\ga_i C_0$ be any points such that
$d(x_i,y_i)\leq c+1$. Up to replacing $\ga_i$ by another
representative of $[\ga_i]$, we may assume that $d(x_i,x_*)\leq D$ and
$d(y_i,\ga_i x_*)\leq D$. Hence $d(x_*,\ga_i x_*)\leq 2D+c+1$ for
every $i$, which contradicts the discreteness of $\Ga$.  \cqfd

\bprop \label{prop:countingdoublecosets}%
Assume that ${\rm Card}\;\Ga x\cap B(x,n)\asymp e^{\delta n}$ for some
(hence every) $x$ in $X$. Then there exists $N$ in $\NN-\{0\}$ such
that
$$
{\rm Card}\;\{r\in R_0\;:\;n \leq D(r)< n+N\}\asymp  
e^{\delta n}\;.
$$
\eprop

\dem As $\delta_0<\delta$, the proof is the same as the proof of
\cite[Theo.~3.4]{HP2}, up to replacing the horoball $H\!B_0$ by $C_0$.
\cqfd

\medskip %
Define $X_0=\Ga_0\backslash X$, and $\partial_\infty X_0=
\Ga_0\backslash(\partial_\infty X-\Lambda\Ga_0)$.  Since
$\Ga_0\backslash C_0$ is compact, and since the closest point map is a
continuous $\Ga_0$-equivariant map from $\partial_\infty
X-\Lambda\Ga_0$ to $C_0$, the space $\partial_\infty X_0$ is compact.
The distance-like map $d_{C_0}$ on $\partial_\infty X-\Lambda\Ga_0$ is
invariant under $\Ga_0$, and we denote by $d_0$ the quotient
distance-like map on $\partial_\infty X_0$, i.e.
$$
d_0(\overline{x},\overline{y})=\inf_{x\in\overline{x},\;
y\in\overline{y}}\;d_{C_0}(x,y)\;.
$$

Let $r=[\ga]$ be an element in $R_0$. Define $\Lambda_r$ (which does
not depend on the representative $\ga$ of $r$) as the image of
$\ga\Lambda\Ga_0$ by the canonical projection $\partial_\infty
X-\Lambda\Ga_0\ra \partial_\infty X_0$.  By Proposition
\ref{prop:equivmalnormal} (2), it follows that $(\Lambda_r)_{r\in
  R_0}$ is a family of pairwise disjoint compact subsets of
$\partial_\infty X_0$.  For every $\epsilon>0$, define
$\N_r(\epsilon)$ as the $\epsilon$-neigbourhood of $\Lambda_r$ in
$\partial_\infty X_0$ for the distance-like map $d_0$. Note that
$\N_r(\epsilon)\subset \N_r(\epsilon')$ if $\epsilon <\epsilon'$.

Let $(\mu_x)_{x\in X}$ be a conformal density of dimension $\delta$
for $\Ga$. Fix a base point $x_0$ in $C_0$. Define
$$
\wt{\mu}_{\Ga_0x_0}= \sum_{\alpha\in\Ga_0}\;\mu_{\alpha x_0}\;.
$$

\blemm\label{lem:abscontmes} %
The map $\wt{\mu}_{\Ga_0x_0}$ is a locally finite Borel measure on
$\partial_\infty X-\Lambda\Ga_0$, which is invariant under $\Ga_0$,
and absolutely continuous with respect to the restriction to
$\partial_\infty X-\Lambda\Ga_0$ of $\mu_x$ for every $x$ in $X$.
\elemm

We denote by $\mu_{\Ga_0x_0}$ the finite Borel measure on the compact
quotient $\partial_\infty X_0$ of $\partial_\infty X-\Lambda\Ga_0$
defined by $\wt{\mu}_{\Ga_0x_0}$.

\medskip \dem Denote by $s\mapsto P_{x_0,\Ga_0}(s)=
\sum_{\alpha\in\Ga_0}\; \;e^{-s d(\alpha x_0,x_0)}$ the Poincar\'e
series of $\Ga_0$ with base point $x_0$.

\noindent
\begin{minipage}{9.6cm}
  ~~~ Let $\xi$ be in $\partial_\infty X - \Lambda\Ga_0$ and $\alpha$
  be in $\Ga_0$. The point $\alpha x_0$ belongs to $C_0$. Hence the
  horosphere centered at $\xi$ passing through $\pi_{C_0}(\xi)$ meets
  the geodesic ray from $\alpha x_0$ to $\xi$ in a point $u$.  As
  $C_0$ is convex and $\pi_{C_0}(\xi)$ is the closed point in $C_0$ to
  $\xi$, by an easy ${\rm CAT}(-1)$ comparison argument, the distance
  $d(u,\pi_{C_0}(\xi))$ is at most $1$ (and even
  $\log\frac{3+\sqrt{5}}{2}\,$). By the triangle inequality,

\end{minipage}
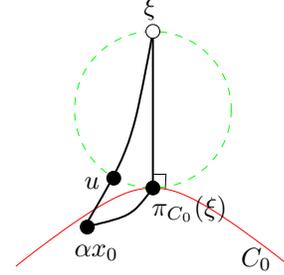
\begin{minipage}{5.3cm}
\begin{center}
\begin{picture}(0,0)%
\includegraphics{fig_plaissaar.pstex}%
\end{picture}%
\setlength{\unitlength}{3236sp}%
\begingroup\makeatletter\ifx\SetFigFont\undefined%
\gdef\SetFigFont#1#2#3#4#5{%
  \reset@font\fontsize{#1}{#2pt}%
  \fontfamily{#3}\fontseries{#4}\fontshape{#5}%
  \selectfont}%
\fi\endgroup%
\begin{picture}(2229,2117)(889,-2069)
\put(1426,-1411){\makebox(0,0)[lb]{\smash{{\SetFigFont{10}{12.0}{\rmdefault}{\mddefault}{\updefault}{\color[rgb]{0,0,0}$u$}%
}}}}
\put(1351,-1936){\makebox(0,0)[lb]{\smash{{\SetFigFont{10}{12.0}{\rmdefault}{\mddefault}{\updefault}{\color[rgb]{0,0,0}$\alpha x_0$}%
}}}}
\put(1944,-1621){\makebox(0,0)[lb]{\smash{{\SetFigFont{10}{12.0}{\rmdefault}{\mddefault}{\updefault}{\color[rgb]{0,0,0}$\pi_{C_0}(\xi)$}%
}}}}
\put(2626,-2011){\makebox(0,0)[lb]{\smash{{\SetFigFont{10}{12.0}{\rmdefault}{\mddefault}{\updefault}{\color[rgb]{0,0,0}$C_0$}%
}}}}
\put(1884,-84){\makebox(0,0)[lb]{\smash{{\SetFigFont{10}{12.0}{\rmdefault}{\mddefault}{\updefault}{\color[rgb]{0,0,0}$\xi$}%
}}}}
\end{picture}%

\end{center}
\end{minipage}

\medskip 
$$
\beta_\xi(\alpha x_0,\pi_{C_0}(\xi))=d(\alpha x_0,u)\geq 
d(\alpha x_0,\pi_{C_0}(\xi))- d(u,\pi_{C_0}(\xi))\;.
$$
Therefore
\begin{align*}
\beta_\xi(\alpha x_0,x_0)&=\beta_\xi(\alpha x_0,\pi_{C_0}(\xi))-
\beta_\xi(x_0,\pi_{C_0}(\xi))\geq d(\alpha x_0,\pi_{C_0}(\xi))-1-
d(x_0,\pi_{C_0}(\xi))\\ &\geq d(\alpha
x_0,x_0)-1-2d(x_0,\pi_{C_0}(\xi))\;,
\end{align*}
where the last inequation is again obtained by the triangle
inequality. Hence
$$
\sum_{\alpha\in\Ga_0}\; \frac{d\mu_{\alpha x_0}}{d\mu_{x_0}}(\xi)=
\sum_{\alpha\in\Ga_0}\;e^{-\delta \beta_\xi(\alpha x_0,x_0)}\leq
e^{1+2d(x_0,\pi_{C_0}(\xi))}\;P_{x_0,\Ga_0}(\delta)\;.
$$
The right hand side, as $\delta>\delta_0$, is a positive continous map
of $\xi\in\partial_\infty X-\Lambda\Ga_0$.  Hence
$\wt{\mu}_{\Ga_0x_0}$ is a locally finite Borel measure on
$\partial_\infty X-\Lambda\Ga_0$.  It is clearly invariant under
$\Ga_0$ by construction and the equivariance property of
$(\mu_x)_{x\in X}$.  As $ \beta_\xi(\alpha x_0,x_0)\leq d(\alpha
x_0,x_0)$, we have, for every $\xi$ in $\partial_\infty
X-\Lambda\Ga_0$,
\begin{equation}\label{eq:derivradnyk}
P_{x_0,\Ga_0}(\delta)\leq
\frac{d\wt{\mu}_{\Ga_0x_0}}{d\mu_{x_0}}(\xi)
\leq e^{1+2d(x_0,\pi_{C_0}(\xi))}\;P_{x_0,\Ga_0}(\delta)\;,
 \end{equation}
hence $\wt{\mu}_{\Ga_0x_0}$ and $\mu_{x_0}$ have the same measure
class on $\partial_\infty X-\Lambda\Ga_0$. \cqfd

\btheo \label{theo:estimatvolepsilneigh}%
There exist constants $c,c'>0$ such that, for every $r$ in $R_0$ and
$\epsilon$ in $]\,0,c'e^{-D(r)}]$,
$$
\frac{1}{c}\;e^{-\delta_0 D(r)}\epsilon^{\delta-\delta_0}\leq
\mu_{\Ga_0x_0}(\N_r(\epsilon))\leq 
c\;e^{-\delta_0D(r)}\epsilon^{\delta-\delta_0}\;.
$$
\etheo

\dem  %
For every double coset $r$ in $R_0$, choose a representative
$\ga_r$ of $r$ such that
$$
d(x_0,\ga_r x_0)=\min_{\alpha,\alpha'\in\Ga_0}
d(x_0,\alpha\ga_r\alpha'x_0)\;.
$$

Denote by $\N_{\epsilon',d'}(A)$ the (closed)
$\epsilon'$-neighbourhood of a subset $A$ for a distance or a
distance-like map $d'$. The subset $\N_{\epsilon,d_{C_0}} (\ga_r
\Lambda\Gamma_0)$ of $\partial_\infty X$ is compact. By Lemma 2.3, it
is contained in $\partial_\infty X-\partial_\infty C_0$, and its
diameter for the distance-like map $d_{C_0}$ tends to $0$ as $D(r)$
tends to $+\infty$ and $\epsilon$ tends to $0$.  Recall that $\Ga_0$
acts isometrically and properly on $\partial_\infty X-
\partial_\infty C_0$ for the distance-like map $d_{C_0}$. 
Hence there exists $N'\in\NN$ and
$c'_1>0$ such that for every $\epsilon$ in $]0,c'_1]$, for every $r$
in $R_0$, we have
$$
{\rm Card}\{\alpha\in\Ga_0\;:\;
\alpha \;\N_{\epsilon,d_{C_0}}(\ga_r\Lambda\Gamma_0)\cap 
\N_{\epsilon,d_{C_0}}(\ga_r\Lambda\Gamma_0)\neq \emptyset\}\leq N'\;.
$$ 

By the construction of $ \mu_{\Ga_0x_0}$, we have, for every $r$ in
$R_0$ and $\epsilon$ in $]0,c'_1]$,
\begin{equation}\label{eq:un}
\frac{1}{N'}\;
\wt{\mu}_{\Ga_0x_0}(\N_{\epsilon,d_{C_0}}(\ga_r\Lambda\Gamma_0))\leq
\mu_{\Ga_0x_0}(\N_r(\epsilon))\leq
\wt{\mu}_{\Ga_0x_0}(\N_{\epsilon,d_{C_0}}(\ga_r\Lambda\Gamma_0))\;.
\end{equation}

As $\Ga_0\backslash C_0$ is compact and by the definition of the
representatives $\ga_r$, there exists $c'_2>0$ such that, for every
$r$ in $R_0$, for every $x\in \ga_r(C_0\cup\partial_\infty C_0)$, the
closest point to $x$ on $C_0$ is at distance at most $c'_2$ from $x_0$
(see also \cite[Lem.~3.5]{HP2}).

Hence, there exists a compact subset $K$ of $\partial_\infty X-
\partial_\infty C_0$ which contains $\ga_r\Lambda\Gamma_0$ for every
$r$ in $R_0$. By Lemma \ref{lem:propdistconv} (2), there exists a
compact subset $K'$ of $\partial_\infty X-
\partial_\infty C_0$ which contains $\N_{\epsilon,d_{C_0}}
(\ga_r\Lambda\Gamma_0)$ for every $\epsilon$ in $]0,c'_1]$ and every
$r$ in $R_0$. Hence by Lemma \ref{lem:propdistconv} (1), there 
exist two constants $c_3^\pm>0$
such that for every $r$ in $R_0$ and $\epsilon\in \;]0,c'_1]$,
\begin{equation}\label{eq:deu}
\N_{c_3^-\epsilon,d_{x_0}}(\ga_r\Lambda\Gamma_0)\subset
\N_{\epsilon,d_{C_0}}(\ga_r\Lambda\Gamma_0)\subset 
\N_{c_3^+\epsilon,d_{x_0}}(\ga_r\Lambda\Gamma_0)\;.
\end{equation}

As $K'$ and $\partial_\infty C_0$ are compact and disjoint, if $c'_1$
is small enough, then there exists a compact subset $K''$ of
$\partial_\infty X- \partial_\infty C_0$ containing
$\N_{c_3^+\epsilon,d_{x_0}} (\ga_r\Lambda\Gamma_0)$ (and hence
$\N_{c_3^-\epsilon,d_{x_0}}(\ga_r\Lambda\Gamma_0)$) for every $r$ in
$R_0$ and $\epsilon$ in $]0,c'_1]$.  By the continuity of $\pi_{C_0}$,
there exists a constant $c'_4>0$ such that for every $r$ in $R_0$ and
$\epsilon$ in $]0,c'_1]$, the subset
$\pi_{C_0}(\N_{c^\pm_3\epsilon,d_{C_0}} (\ga_r\Lambda\Gamma_0))$ is
contained in the ball of center $x_0$ and radius $c'_4$.

By the definition of the representatives $\ga_r$, for every $r$ in
$R_0$, for every $\xi\in \ga_r\partial_\infty C_0$, the point $\ga_r
x_0$ is at distance at most a constant from the geodesic between $x_0$
and $\xi$ (see also \cite[Lem.~3.5]{HP2}). Recall that for every
$\eta,\eta'$ in $\partial_\infty X$, if $d_{x_0}(\eta,\eta')\leq
\epsilon'$, then the geodesic rays $[x_0,\eta[$ and $[x_0,\eta'[$
remain at distance bounded by a universal constant at least during a
time $-\log \epsilon'$.  Hence, if $c'\leq c'_1$ is small enough and
$\epsilon\leq c'e^{-D(r)}$, then every geodesic ray from $x_0$ to a
point $\xi$ in $\N_{c_3^\pm\epsilon,d_{x_0}} (\ga_r\Lambda\Gamma_0)$
passes at distance less than a constant from $\ga_rx_0$.  This has two
consequences.

$\bullet$ First, using the change of base point formula for the visual
distances, there exist two constants $c_5^\pm>0$ such that for every
$\epsilon\leq c'\,e^{-D(r)}$,
\begin{equation}\label{eq:tro}
\N_{c_3^+\epsilon,d_{x_0}}(\ga_r\Lambda\Gamma_0) \subset
\N_{c_5^+e^{D(r)}\epsilon,d_{\ga_r x_0}}(\ga_r\Lambda\Gamma_0)
\;\;{\rm and}\;\;
\N_{c_5^-e^{D(r)}\epsilon,d_{\ga_r x_0}}(\ga_r\Lambda\Gamma_0)
\subset \N_{c_3^-\epsilon,d_{x_0}}(\ga_r\Lambda\Gamma_0)\;.
\end{equation}

$\bullet$ Second, for every $\xi$ in $\N_{c_3^\pm\epsilon,d_{x_0}}
(\ga_r\Lambda\Gamma_0)$, the number $|\beta_\xi(x_0,\ga_r
x_0)-d(x_0,\ga_r x_0)|$ is bounded by a constant. Hence there exist
constants $c_{6}^\pm>0$ such that for every $r$ in $R_0$ and for every
$\xi$ in $\N_{c_3^\pm\epsilon,d_{x_0}} (\ga_r\Lambda\Gamma_0)$,
\begin{equation}\label{eq:qua}
c_{6}^-\;e^{-\delta D(r)}\leq 
\frac{d\mu_{x_0}}{d\mu_{\ga_r x_0}}(\xi)
\leq c_{6}^+\;e^{-\delta D(r)}\;.
\end{equation}

By the Radon-Nykodim derivative estimates in Equation
\eqref{eq:derivradnyk} and the definition of $c'_4$, there exist
constants $c_7^\pm>0$ such that for every $\epsilon$ in $]0,c'_1]$,
every $r$ in $R_0$, and every $\xi$ in $\N_{c_3^\pm \epsilon,d_{x_0}}
(\ga_r\Lambda\Gamma_0)$,
\begin{equation}\label{eq:cin}
c_7^-\leq \frac{d\wt{\mu}_{\Ga_0x_0}}{d\mu_{x_0}}(\xi)\leq c_7^+\;.
\end{equation}

By Sullivan's shadow lemma (see for instance \cite[Lem.~1.3]{Rob}),
for every constant $c_8'>0$ big enough, there exist constants
$c_9^\pm>0$ such that, for every $\ga$ in $\Ga$,
\begin{equation}\label{eq:six}
c_9^-e^{-\delta d(x_0,\ga x_0)}\leq 
\mu_{x_0}(\O_{x_0}B(\ga x_0,c_8'))
\leq c_9^+e^{-\delta d(x_0,\ga x_0)}\;.
\end{equation}

For every $t\geq 0$, define $\Ga_0[t]=\{\alpha\in\Ga_0\;:\;
d(x_0,\alpha x_0)\leq t\}$. For every $\epsilon'\in\;]\,0,1]$ and
$\kappa>0$, define
$$
A^+_{\epsilon',\kappa}= 
\Ga_0[-\log \epsilon'+\kappa]-\Ga_0[-\log \epsilon' -\kappa]
\;\;{\rm and}\;\;  A^-_{\epsilon',\kappa}= 
\Ga_0[-\log \epsilon' + 2 \kappa]-\Ga_0[-\log \epsilon' + \kappa]
\;.
$$ 
Let $\epsilon'\in\;]\,0,1]$, $\eta\in\partial_\infty X$ and
$\eta'\in\Lambda\Ga_0$ be such that $d_{x_0}(\eta,\eta')\leq
\epsilon'$.  Let $u$ be the point of $[x_0,\eta[$ at distance $-\log
\epsilon'$ from $x_0$. By the definition of $d_{x_0}$ and the
properties of the geodesic rays in a ${\rm CAT}(-1)$ metric space,
there exists a universal constant $c''_8$ such that $\eta'$ belongs to
$\O_{x_0} B(u,c''_8)$. Let $c^+_{10}>0$ be at least the (finite)
diameter of $\Ga_0\backslash C_0$. Since $\partial_\infty
C_0=\Lambda\Ga_0$ and by convexity, the geodesic ray $[x_0,\eta[$ is
contained in $C_0$. Hence there exists $\alpha$ in $\Ga_0$ such that
$d(u,\alpha x_0)\leq c^+_{10}$. Let $c'_8$ be big enough (at least
$c''_8+c^+_{10}$ and such that Equation \eqref{eq:six} holds).  Then,
by the triangle inequality, $B(u,c''_8)$ is contained in $B(\alpha
x_0,c'_8)$. Note that $-\log \epsilon'-c^+_{10}\leq d(x_0,\alpha
x_0)\leq -\log \epsilon'+c^+_{10} $.  Therefore, for every $\epsilon'$
in $]0,1]$, we have
\begin{equation}\label{eq:sep}
\N_{\epsilon',d_{x_0}} (\Lambda\Gamma_0)\subset 
\bigcup_{\alpha\in A^+_{\epsilon',c^+_{10}}} 
\O_{x_0}B(\alpha x_0,c_8')\;.
\end{equation}
As $\Ga_0$ is convex-cocompact, there exists a constant $c'''_8>0$
such that for every $\alpha$ in $\Gamma_0$, the segment $[x_0,\alpha
x_0]$ is at distance at most $c'''_8$ from a geodesic ray starting
from $x_0$ and contained in $C_0$. Assume that $c^-_{10}>0$ is at
least $c'''_8+c'_8$.  Let $\epsilon'\in\;]\,0,1]$, $\alpha\in
A^-_{\epsilon',c^-_{10}}$ and $\eta'\in \O_{x_0}B(\alpha x_0,c_8')$.
Let $v$ be a point on $[x_0,\eta'[$ at distance at most $c'_8$ from
$\alpha x_0$. Let $\eta\in \partial_\infty C_0$ and $u\in [x_0,\eta[$
be such that $d(u,\alpha x_0)\leq c'''_8$, which exist by the
definition of $c'''_8$.  Then by the definition of $d_{x_0}$ and the
triangle inequality, we have
$$
d_{x_0}(\eta,\eta')\leq e^{\frac{1}{2}(d(u,v)-d(x_0,u)-d(x_0,v))}\leq
e^{c'_8+c'''_8-d(x_0,\alpha x_0)} \leq \epsilon'\;,
$$
since $\alpha\in A^-_{\epsilon',c^-_{10}}$. Therefore, for
every $\epsilon'$ in $]0,1]$, we have
$$
\bigcup_{\alpha\in A^-_{\epsilon',c_{10}^-}} 
\O_{x_0}B(\alpha x_0,c_8')\subset 
\N_{\epsilon',d_{x_0}} (\Lambda\Gamma_0)\;.
$$
If $c_{10}^+$ and then $c_{10}^-$ are big enough, as $\Ga_0$ is
convex-cocompact
(by for instance \cite{Rob} if $\Ga_0$ is non elementary, and even if
$\delta_0=0$, since then, by the assumptions, $\Ga_0$ contains an
hyperbolic element generating a finite index (infinite cyclic)
subgroup), note that there exist constants $c_{11}^{\pm}>0$ such that
for every $\epsilon'$ in $]0,1]$, we have
\begin{equation}\label{eq:neu}
{\rm Card~} A^-_{\epsilon',c_{10}^-}\geq c_{11}^-\;{(\epsilon')}^{-\delta_0}
\;\;{\rm and}\;\;
{\rm  Card~}A^+_{\epsilon',c_{10}^+} 
\leq c_{11}^+\;{(\epsilon')}^{-\delta_0}\;.
\end{equation}

Let $A^*_{\epsilon'}$ be a maximal subset of $A^-_{\epsilon',c_{10}^-}$
such that the shadows $\O_{x_0}B(\alpha x_0,c_8')$ for $\alpha$ in
$A^*_{\epsilon'}$ are pairwise disjoint.  By maximality, for every
$\alpha$ in $A^-_{\epsilon',c_{10}^-}$, there exists $\alpha'$ in
$A^*_{\epsilon'}$ such that $\alpha x_0$ and $\alpha' x_0$ are at
bounded distance. Hence there exists a constant $c_{12}'>0$ such that
${\rm Card~}A^*_{\epsilon'}\geq c_{12}'{\rm
~Card~}A^-_{\epsilon',c_{10}^-}$.

Let us now prove the upper bound in Theorem
\ref{theo:estimatvolepsilneigh}.  Let $r$ in $R_0$ and $\epsilon$ in
$]\,0,c'\,e^{-D(r)}]$. Note that if $c'>0$ is small enough, then 
$\epsilon\leq c'_1$ and
$c_5^\pm \,e^{D(r)}\epsilon\leq c_5^\pm\,c'\leq 1$. We have
\begin{align*}
\mu_{\Ga_0x_0}(\N_r(\epsilon))
& \leq
\wt{\mu}_{\Ga_0x_0}(\N_{\epsilon,d_{C_0}}(\ga_r\Lambda\Gamma_0)) 
\leq
\wt{\mu}_{\Ga_0x_0}( \N_{c_3^+\epsilon,d_{x_0}}(\ga_r\Lambda\Gamma_0)) 
& {\rm by}\; \eqref{eq:un}\;{\rm and}\; \eqref{eq:deu}\\ &\leq
c_7^+\mu_{x_0}(\N_{c_3^+\epsilon,d_{x_0}}(\ga_r\Lambda\Gamma_0)) 
&{\rm by}\; \eqref{eq:cin}\\ & \leq 
c_7^+c_{6}^+e^{-\delta D(r)} \mu_{\ga_r x_0}(
\N_{c_5^+e^{D(r)}\epsilon,d_{\ga_r x_0}}(\ga_r\Lambda\Gamma_0))
&{\rm by}\; \eqref{eq:qua}\;{\rm and}\; \eqref{eq:tro}\\ &=
c_7^+c_{6}^+e^{-\delta D(r)}
\mu_{x_0}(\N_{c_5^+e^{D(r)}\epsilon,d_{x_0}}(\Lambda\Gamma_0))
&{\rm by~invariance}\\ &
\leq c_7^+c_{6}^+e^{-\delta D(r)}\sum_{\alpha\in
A^+_{c_5^+e^{D(r)}\epsilon,c_{10}^+}} 
\mu_{x_0}(\O_{x_0}B(\alpha x_0,c_8'))
&{\rm by}\; \eqref{eq:sep}\\&
\leq c_7^+c_{6}^+ e^{-\delta D(r)} c_9^+
e^{-\delta(-\log(c_5^+e^{D(r)}\epsilon) -c_{10}^+)} c_{11}^+
  (c_5^+e^{D(r)}\epsilon)^{-\delta_0}
&{\rm by}\; \eqref{eq:six}\;{\rm and}\; \eqref{eq:neu}\\&=
c_{13}^+ \;e^{-\delta_0 D(r)}  \epsilon^{\delta-\delta_0}\;,
\end{align*}
for some constant $c_{13}^+>0$, which proves the upper bound.

Similarly for the lower bound, 
\begin{align*}
\mu_{\Ga_0x_0}(\N_r(\epsilon))& \geq\frac{1}{N'}\;
\wt{\mu}_{\Ga_0x_0}( \N_{\epsilon,d_{C_0}}(\ga_r\Lambda\Gamma_0)) 
\geq
\frac{1}{N'}\;\wt{\mu}_{\Ga_0x_0}
( \N_{c_3^-\epsilon,d_{x_0}}(\ga_r\Lambda\Gamma_0)) 
\\& \geq 
\frac{c_7^-c_{6}^-}{N'}\;e^{-\delta D(r)}
\mu_{x_0}(\N_{c_5^-e^{D(r)}\epsilon,d_{x_0}}(\Lambda\Gamma_0))
\\&\geq
 \frac{c_7^-c_{6}^-}{N'}\;e^{-\delta D(r)}\mu_{x_0}\Big(\bigcup_{\alpha\in
    A^*_{c_5^-e^{D(r)}\epsilon}} \O_{x_0}B(\alpha x_0,c_8')\Big)
\\&=
\frac{c_7^-c_{6}^-}{N'}\;e^{-\delta D(r)}\sum_{\alpha\in
  A^*_{c_5^-e^{D(r)}\epsilon}} \mu_{x_0}(\O_{x_0}B(\alpha x_0,c_8'))
\\&
\geq \frac{c_7^-c_{6}^-}{N'}\; e^{-\delta D(r)} \;c_9^-\;
e^{-\delta(-\log(c_5^-e^{D(r)}\epsilon) +2c_{10}^-)} \;c_{12}'\;c_{11}^-
  (c_5^-e^{D(r)}\epsilon)^{-\delta_0}\\&
= c_{13}^- \;e^{-\delta_0 D(r)}  \epsilon^{\delta-\delta_0}\;,
\end{align*}
for some constant $c_{13}^->0$, which proves the result.  
\cqfd

\blemm \label{lem:verifsixoldthree} %
For every $N\in\NN-\{0\}$, there exists $c''>0$ such that for every
$n$ in $\NN$, for every distinct $r$ and $r'$ in $R_0$ such that
$D(r)$ and $D(r')$ belong to $[nN,(n+1)N[\,$, the subsets 
$\N_r(c''\,e^{- n N})$ and $\N_{r'}(c''\,e^{- n N})$ are disjoint.  
\elemm

\dem 
Let $N\in\NN-\{0\}$, and $c''\leq 1$ be small enough, to be determined
during the proof. Assume by absurd that there exists $n$ in $\NN$,
distinct $r$ and $r'$ in $R_0$ such that $D(r),D(r')\in
[nN,(n+1)N[\,$, and that the subsets $\N_r(c''\,e^{- n N})$ and
$\N_{r'}(c''\,e^{- n N})$ have non empty intersection.  Then, there
exist representatives $\ga,\ga'$ of the double cosets $r,r'$ and
points $\xi,\xi'$ in $\ga\Lambda\Ga_0, \ga' \Lambda\Ga_0$
respectively, and an element $\eta$ in $\partial_\infty
X-\partial_\infty C_0$ which is different from $\xi,\xi'$, such that
$d_{C_0}(\xi,\eta)$ and $d_{C_0}(\xi',\eta)$ are at most $c''\,e^{- n
  N}$, and in particular at most $c"$.

Since there are only finitely many $r$'s with $D(r)$ less than a
constant, and since the subsets $\ga\Lambda\Ga_0$ for $\ga$ in
$(\Ga-\Ga_0)/\Ga_0$ are pairwise disjoint (by Proposition
\ref{prop:equivmalnormal} (2)) closed subsets, we may assume that
$D(r)$ and $D(r')$ are bigger than any given constant $c''_1>0$. In
particular, $D(r)$ and $D(r')$ are positive.

By Lemma \ref{lem:propdistconv} (4), there exists a universal constant
$c''_2\geq 1$ such that if $d_{C_0}(\eta',\eta'')\leq 1/c''_2$, then
the geodesic line between $\eta'$ and $\eta''$ is disjoint from $C_0$,
and the length of the common perpendicular segment between
$]\eta',\eta''[$ and $C_0$ is at most $-\log d_{C_0}(\eta',\eta'')
+c''_2$ and at least $-\log d_{C_0}(\eta',\eta'') -c''_2$.  Assume
that $c''\leq 1/c''_2$.

\medskip
\noindent
\begin{minipage}{7.5cm}
  ~~~Let $p_\xi,p_{\xi'},p_\eta$ be the closest point on $C_0$ to
  $\xi,\xi',\eta$ respectively. Let $[x,y]$ (resp.~$[x',y']$;
  $[x_\xi,y_\xi]$; $[x_{\xi'},y_{\xi'}]$) be the common perpendicular
  between $C_0$ and $\ga C_0$ (resp.~$C_0$ and $\ga' C_0$; $C_0$ and
  $]\xi,\eta[$; $C_0$ and $]\xi',\eta[$), with $x,x',x_\xi,x_{\xi'}$
  in $C_0$.  Let $z,z'$ be the closest point to $y,y'$ on $[p_\xi,\xi[
  , [p_{\xi'},\xi'[$ respectively. Let $v,u,v',u'$ be the closest
  points to $y_\xi,y_\xi,y_{\xi'},y_{\xi'}$ on $[p_\xi,\xi[
  ,[p_\eta,\eta[, [p_{\xi'},\xi'[,[p_\eta,\eta[$ respectively (see the
  picture on the right).  Let $z_\xi,v_\ga,w_{\ga'}$ be the closest
  point to $z,v,u$ on $[x_\xi,y_\xi],\ga C_0,\ga'C_0$ respectively.
  We have $d(x,y)=D(r)$, $d(x',y')=D(r')$.
\end{minipage}
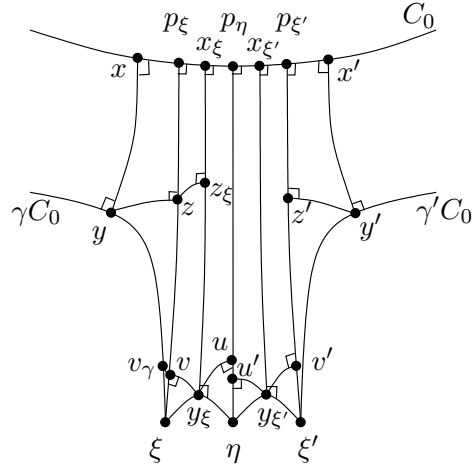
\begin{minipage}{6.3cm}
\begin{center}
\begin{picture}(0,0)%
\includegraphics{fig_separprojnomb.pstex}%
\end{picture}%
\setlength{\unitlength}{3729sp}%
\begingroup\makeatletter\ifx\SetFigFont\undefined%
\gdef\SetFigFont#1#2#3#4#5{%
  \reset@font\fontsize{#1}{#2pt}%
  \fontfamily{#3}\fontseries{#4}\fontshape{#5}%
  \selectfont}%
\fi\endgroup%
\begin{picture}(3093,3139)(1420,-7472)
\put(4276,-4516){\makebox(0,0)[lb]{\smash{{\SetFigFont{11}{13.2}{\rmdefault}{\mddefault}{\updefault}{\color[rgb]{0,0,0}$C_0$}%
}}}}
\put(2212,-5950){\makebox(0,0)[lb]{\smash{{\SetFigFont{11}{13.2}{\rmdefault}{\mddefault}{\updefault}{\color[rgb]{0,0,0}$y$}%
}}}}
\put(3671,-6820){\makebox(0,0)[lb]{\smash{{\SetFigFont{11}{13.2}{\rmdefault}{\mddefault}{\updefault}{\color[rgb]{0,0,0}$v'$}%
}}}}
\put(2842,-7152){\makebox(0,0)[lb]{\smash{{\SetFigFont{11}{13.2}{\rmdefault}{\mddefault}{\updefault}{\color[rgb]{0,0,0}$y_\xi$}%
}}}}
\put(3843,-4900){\makebox(0,0)[lb]{\smash{{\SetFigFont{11}{13.2}{\rmdefault}{\mddefault}{\updefault}{\color[rgb]{0,0,0}$x'$}%
}}}}
\put(3240,-4698){\makebox(0,0)[lb]{\smash{{\SetFigFont{11}{13.2}{\rmdefault}{\mddefault}{\updefault}{\color[rgb]{0,0,0}$x_{\xi'}$}%
}}}}
\put(2903,-4690){\makebox(0,0)[lb]{\smash{{\SetFigFont{11}{13.2}{\rmdefault}{\mddefault}{\updefault}{\color[rgb]{0,0,0}$x_\xi$}%
}}}}
\put(3991,-5927){\makebox(0,0)[lb]{\smash{{\SetFigFont{11}{13.2}{\rmdefault}{\mddefault}{\updefault}{\color[rgb]{0,0,0}$y'$}%
}}}}
\put(2801,-5792){\makebox(0,0)[lb]{\smash{{\SetFigFont{11}{13.2}{\rmdefault}{\mddefault}{\updefault}{\color[rgb]{0,0,0}$z$}%
}}}}
\put(3533,-5845){\makebox(0,0)[lb]{\smash{{\SetFigFont{11}{13.2}{\rmdefault}{\mddefault}{\updefault}{\color[rgb]{0,0,0}$z'$}%
}}}}
\put(3173,-6863){\makebox(0,0)[lb]{\smash{{\SetFigFont{11}{13.2}{\rmdefault}{\mddefault}{\updefault}{\color[rgb]{0,0,0}$u'$}%
}}}}
\put(3008,-6687){\makebox(0,0)[lb]{\smash{{\SetFigFont{11}{13.2}{\rmdefault}{\mddefault}{\updefault}{\color[rgb]{0,0,0}$u$}%
}}}}
\put(3311,-7160){\makebox(0,0)[lb]{\smash{{\SetFigFont{11}{13.2}{\rmdefault}{\mddefault}{\updefault}{\color[rgb]{0,0,0}$y_{\xi'}$}%
}}}}
\put(3101,-7400){\makebox(0,0)[lb]{\smash{{\SetFigFont{11}{13.2}{\rmdefault}{\mddefault}{\updefault}{\color[rgb]{0,0,0}$\eta$}%
}}}}
\put(2594,-7399){\makebox(0,0)[lb]{\smash{{\SetFigFont{11}{13.2}{\rmdefault}{\mddefault}{\updefault}{\color[rgb]{0,0,0}$\xi$}%
}}}}
\put(3580,-7403){\makebox(0,0)[lb]{\smash{{\SetFigFont{11}{13.2}{\rmdefault}{\mddefault}{\updefault}{\color[rgb]{0,0,0}$\xi'$}%
}}}}
\put(2992,-5665){\makebox(0,0)[lb]{\smash{{\SetFigFont{11}{13.2}{\rmdefault}{\mddefault}{\updefault}{\color[rgb]{0,0,0}$z_\xi$}%
}}}}
\put(2449,-6827){\makebox(0,0)[lb]{\smash{{\SetFigFont{11}{13.2}{\rmdefault}{\mddefault}{\updefault}{\color[rgb]{0,0,0}$v_\ga$}%
}}}}
\put(2771,-6823){\makebox(0,0)[lb]{\smash{{\SetFigFont{11}{13.2}{\rmdefault}{\mddefault}{\updefault}{\color[rgb]{0,0,0}$v$}%
}}}}
\put(1435,-5633){\makebox(0,0)[lb]{\smash{{\SetFigFont{11}{13.2}{\rmdefault}{\mddefault}{\updefault}{\color[rgb]{0,0,0}${}$}%
}}}}
\put(2691,-4533){\makebox(0,0)[lb]{\smash{{\SetFigFont{11}{13.2}{\rmdefault}{\mddefault}{\updefault}{\color[rgb]{0,0,0}$p_\xi$}%
}}}}
\put(3070,-4536){\makebox(0,0)[lb]{\smash{{\SetFigFont{11}{13.2}{\rmdefault}{\mddefault}{\updefault}{\color[rgb]{0,0,0}$p_\eta$}%
}}}}
\put(3439,-4532){\makebox(0,0)[lb]{\smash{{\SetFigFont{11}{13.2}{\rmdefault}{\mddefault}{\updefault}{\color[rgb]{0,0,0}$p_{\xi'}$}%
}}}}
\put(2335,-4871){\makebox(0,0)[lb]{\smash{{\SetFigFont{11}{13.2}{\rmdefault}{\mddefault}{\updefault}{\color[rgb]{0,0,0}$x$}%
}}}}
\put(1678,-5821){\makebox(0,0)[lb]{\smash{{\SetFigFont{11}{13.2}{\rmdefault}{\mddefault}{\updefault}{\color[rgb]{0,0,0}$\ga  C_0$}%
}}}}
\put(4368,-5818){\makebox(0,0)[lb]{\smash{{\SetFigFont{11}{13.2}{\rmdefault}{\mddefault}{\updefault}{\color[rgb]{0,0,0}$\ga'C_0$}%
}}}}
\end{picture}%

\end{center}
\end{minipage}
\medskip

Assume that $c''\leq e^{-c''_1-c''_2}$, so that $d(x_\xi,y_\xi)$,
$d(x_{\xi'},y_{\xi'})$ are at least $-\log(c'') -c''_2\geq c''_1$. By
the convexity of $C_0$ and quasi-geodesic arguments, if $c''_1$ is
bigger than some universal constant, then there exists a universal
constant $c''_3$ such that the distances $d(z,y) $, $d(x,p_\xi)$, $
d(z',y')$, $ d(x',p_{\xi'}) $, $d(v,y_\xi)$, $ d(u,y_\xi)$,
$d(v',y_{\xi'})$, $ d(u',y_{\xi'})$, $ d(p_\xi,x_\xi)$,
$d(p_\eta,x_\xi)$, $d(p_\eta,x_{\xi'})$, $d(p_{\xi'},x_{\xi'})$ are at
most $c''_3$.  Hence
\begin{align*}
 d(z,v) &  \geq d(z,y_\xi)-d(y_\xi,v)
 \geq d(z_\xi,y_\xi)-c''_3\geq
d(y_\xi,x_\xi)-d(x_\xi,z_\xi)-c''_3\\
&\geq
(-\log(c''e^{-nN})-c''_2)-(d(x,y)+2c''_3)-c''_3\\
& \geq 
-\log c'' +nN-c''_2-N(n+1)-3\,c''_3=-\log c'' -N-c''_2-3c''_3\;.
\end{align*}
In particular, if $c''$ is small, then the points $p_\xi,z,v$ are in
this order on $[p_\xi,\xi[$. So that by convexity $d(y,v_\ga)\geq
d(z,v)-c''_3\geq -\log c'' -N-c''_2-4c''_3$, which is big if $c''$ is
small. Similarly, $p_{\xi'},z',v'$ are in this order on
$[p_{\xi'},\xi'[$.

Up to permuting $\xi,\xi'$, we may assume that $p_\eta,u,u',\eta$ are
in this order on the geodesic ray $[p_\eta,\eta[$.  By convexity
properties of the distance, $d(u,w_{\ga'})\leq 3c''_3$. Hence as
above, there exists a constant $c''_4>0$ such that $d(y',w_{\ga'})\geq
-\log c'' -c''_4$.

Consider the geodesic hexagon with vertices $x,y,v_\ga,w_{\ga'}, y',
x'$. We have $d(x,x')\leq 6c''_3$, $d(v_\ga,w_{\ga'})\leq 6c''_3$, and
$[x,y]$ (resp.~$[x',y']$) is the shortest segment between points of
$[x,x']$ and $[y,v_\ga]$ (resp.~$[y',w_{\ga'}]$). Furthermore $d(x,y),
d(x',y')$ may be assumed to be bigger than any constant and differ by
at most a constant (that is by $N$), and $d(y,v_\ga), d(y',w_{\ga'})$
are bigger than $-\log c''$ minus a constant. For instance by using
techniques of approximation by trees (see for example \cite[page
33]{GH}) on the above hexagon, the geodesic segments $[y,v_\ga]$ and
$[y',w_{\ga'}]$, contained respectively in $\ga C_0$ and $\ga' C_0$,
are arbitrarily long if $c''$ is small enough; moreover their first
endpoints and last endpoints are at bounded distance.  Let $\epsilon
=1$ and $\kappa(\epsilon)$ be given by Proposition
\ref{prop:equivmalnormal} (4). Hence by hyperbolicity, the
$\epsilon$-neighborhoods of $\ga C_0$ and $\ga' C_0$ meet in a segment
of length that can be made bigger than $\kappa(\epsilon)$ if $c''$ is
small enough. This is a contradiction to Proposition
\ref{prop:equivmalnormal} (4).  \cqfd

\medskip %
A map $\psi:[0,+\infty[\;\ra\;]0,+\infty[$ is called {\it slowly
  varying} (see \cite{Sul}) if it is measurable and if there exist
constants $B>0$ and $A\geq 1$ such that for every $x,y$ in $\RR_+$, if
$|x-y|\leq B$, then $\psi(y)\leq A\,\psi(x)$.  Recall (see for
instance \cite[Sec.~5]{HP2}) that this implies that $\psi$ is locally
bounded, hence it is locally integrable; also, if $\log \psi$ is
Lipschitz, then $\psi$ is slowly varying; and for every
$N\in\NN-\{0\}$ and $\epsilon>0$, the series $\sum_{n=0}^{\infty}
\psi(Nn)^\epsilon$ converges if and only if the integral
$\int_{0}^{\infty} \psi(t)^\epsilon dt$ converges.

The following statement is the main technical step towards our
Khintchine-type theorem for the spiraling of geodesic rays in
$\Ga\backslash X$ around $\Ga_0\backslash C_0$. It gives a $0$-$1$
measure result for the approximation of points in the limit set of
$\Ga$ by points of the orbit under $\Ga$ of the limit set of $\Ga_0$.

\btheo \label{theo:maintechnical} %
Let $X,\Ga,\Ga_0,\delta,\delta_0, (\mu_x)_{x\in X}$ be as above.
Assume furthermore, for some (hence any) $x$ in $X$, that $\mu_x$ is
ergodic for the action of $\Ga$, and that ${\rm Card}\;\Ga x\cap
B(x,n)\asymp e^{\delta n}$.  Let $f:[0,+\infty[\;\ra \;]0,+\infty[$ be
a slowly varying map.

If $\int_{1}^{+\infty} f(t)^{\delta-\delta_0}\;dt$ converges
(resp.~diverges), then $\mu_{\Ga_0x_0}$-almost no (resp.~every) point
of $\partial_\infty X_0$ belongs to infinitely many
$\N_r\bigl(f(D(r))e^{-D(r)}\bigr)$ where $r\in R_0$.  
\etheo

\rem In addition to the hypotheses on $X$ in this theorem, assume in
this remark that $X$ is a Riemannian manifold with constant sectional
curvature $-1$, that $\Ga$ is convex-cocompact and that $\Ga_0$ is the
stabilizer of a geodesic line.  Then up to some rewriting, this result
is already known, see for instance \cite{DMPV} or the recent
\cite{BV}. But even in this particular case, our techniques are very
different from the ones of \cite{DMPV,BV}.

\medskip
\dem %
By a similar reduction as in \cite[Lem.~5.2]{HP2}, we may assume
that $f\leq 1$.  Define $g=-\log f:[0,+\infty[\;\ra[0,+\infty[$.

We apply Theorem \ref{theo:resucborelcantelli} with $Z=\partial_\infty
X_0$, $\mu=\mu_{\Ga_0x_0}$, $I=R_0$, and, for every $r$ in $R_0$,
$n\in \NN$ and $\epsilon>0$, with $B_r(\epsilon)= \N_r(\epsilon)$,
$I_n=\{r\in R_0\;:\;Nn\leq D(r)<N(n+1)\}$ where $N$ is as in
Proposition \ref{prop:countingdoublecosets}, and
$n_r=E[\frac{D(r)}{N}]$ where $E$ denotes the integer part. Define,
for every $n$ in $\NN$ and $\epsilon>0$,
$$
f_1(n)=e^{\delta n N},\;\;
f_2(n)=c_2\,e^{- n N},\;\;
f_3(n)=c_2\,e^{-(nN +g(nN))},\;\;
f_4(n)=c_2^{\delta_0-\delta}\,e^{-\delta_0 nN}\;,$$
$$f_5(\epsilon)=\epsilon^{\delta-\delta_0},
$$
where $c_2$ is a small enough positive constant. In particular, we
assume that $c_2$ is less than $c'e^{-N}$, where $c'$ is the constant
defined in Theorem \ref{theo:estimatvolepsilneigh}, and less than the
constant $c''$ defined in Lemma \ref{lem:verifsixoldthree}. Note that
$$
f_1(n)f_4(n)f_5(f_3(n))= 
e^{-(\delta-\delta_0)g(Nn)}= 
f(Nn)^{\delta-\delta_0}\;.
$$  
Hence, as $f$ is slowly varying, the series $\sum_{n\in\NN}
f_1(n)f_4(n)f_5(f_3(n))$ converges if and only if the integral
$\int_{1}^{+\infty} f^{\delta-\delta_0}$ converges.

Note that $B_r(\epsilon)$ is measurable and non-decreasing in
$\epsilon$, and that $I_n$ is finite by Lemma
\ref{lem:depthsdiscrete}.  Assumption (1) of Theorem
\ref{theo:resucborelcantelli} is satisfied since $g$ is non negative.
The assumptions (2) and (3) are easily verified. Assumption (4) follows
from Proposition \ref{prop:countingdoublecosets}. Assumption (5)
follows from Theorem \ref{theo:estimatvolepsilneigh} and the first
assumption on $c_2$. Assumption (6) is satisfied by Lemma
\ref{lem:verifsixoldthree} and the second assumption on $c_2$.

Let us check that Assumption (7) of Theorem
\ref{theo:resucborelcantelli} is also satisfied.

\medskip 
Let $r,r'\in R_0$ with $n=n_r<m=n_{r'}$ such that $\N_r(f_3(n))$ and
$\N_{r'}(f_3(m))$ meet. Hence, there exists representatives $\ga,\ga'$
of $r,r'$, points $\xi,\xi'$ in $\ga\Lambda\Ga_0,\ga'\Lambda\Ga_0$
respectively, and $\eta\in\partial_\infty X-\Lambda\Ga_0$ such that
$d_{C_0}(\xi,\eta)\leq f_3(n)$ and $d_{C_0}(\xi',\eta)\leq f_3(m)$.
Let us prove that there exists a big enough constant $\lambda>0$ such
that $\N_{r'}(f_2(m))$ is contained in $\N_r(\lambda f_3(n))$.

Recall that there are only finitely many $r$'s with $D(r)$ less than a
constant. As $d_{C_0}$ is bounded on $\partial_\infty
X_0\times\partial_\infty X_0$, the $\epsilon$-neighbourhood for
$d_{C_0}$ of any non-empty set covers $\partial_\infty X_0$ if
$\epsilon$ is big enough. Hence we may assume, if $\lambda$ is big
enough, that $D(r)$ and $D(r')$ are bigger than any given constant
$c_6>0$. In particular, $D(r)$ and $D(r')$ are positive.

Let $p_\xi,p_{\xi'},p_\eta, x,y,x',y',x_\xi,y_\xi,x_{\xi'},y_{\xi'},
z,z',v,u,v',u',v_\ga$ be as in the proof of Lemma
\ref{lem:verifsixoldthree} and its picture. Let $v_{\ga'}$ be the
closest point to $v'$ on $\ga'C_0$.

As in the proof of Lemma \ref{lem:verifsixoldthree}, if $c_6$ is
bigger than a universal constant and if $c_2$ is small enough, then
there exists a universal constant $c_7$ such that the following
distances $d(z,y) $, $d(x,p_\xi,)$, $ d(z',y')$, $ d(x',p_\xi') $,
$d(v,y_\xi)$, $ d(u,y_\xi)$, $d(v',y_{\xi'})$, $ d(u',y_{\xi'})$, $
d(p_\xi,x_\xi)$, $d(p_\eta,x_\xi)$, $d(p_\eta,x_{\xi'})$,
$d(p_{\xi'},x_{\xi'})$, $d(v,v_\ga)$, $d(v',v_{\ga'})$ are at most
$c_7$. Furthermore, $p_\xi,z,v,\xi$ are in this order on $[p_\xi,\xi[$
and $p_{\xi'},z',v',\xi'$ are in this order on $[p_{\xi'},\xi'[\,$,
and $d(y,v_\ga), d(y',v_{\ga'})$ may be taken bigger than any given
constant if $c_2$ is small enough.

Say that a point $p$ is {\it above} $q$ (resp.~{\it below $q$ by at
  most some constant $h>0$)} with respect to $C_0$ if $d(p,C_0)\geq
d(q,C_0)$ (resp.~$d(q,C_0)\geq d(p,C_0)\geq d(q,C_0)-h$).  As $m>n$,
the point $y'$ is above $y$ or below $y$ by at most some universal
constant.  If the point $y'$ was below $u$ by more than some big
constant, then, if $c_2$ is small enough, some long subsegment of
$[y',v_{\ga'}]$ would have its endpoints at distance at most a few
$c_7$'s from the endpoints of some subsegment of $[y,v_{\ga}]$, and as
in the end of the proof of Lemma \ref{lem:verifsixoldthree}, this
would contradict Proposition \ref{prop:equivmalnormal} (4). Therefore
the point $y'$ is either above, or below only by a some constant, the
point $u$ and hence $y_{\xi}$.  So that for every $\lambda''>0$, there
exists $\lambda'>0$ such that the shadow (seen from $p_\xi$) of the
ball of center $y_{\xi}$ and radius $\lambda'>0$ contains the shadow
of the ball of center $y'$ and radius $\lambda''>0$. Note that if
$\lambda''$ is big enough, then the shadow of $B(y',\lambda'')$
contains $\N_{r'}(f_2(m))$, as seen in the proof of Lemma
\ref{lem:verifsixoldthree}.  But if $\lambda$ is big enough, then
$\N_r(\lambda f_3(n))$ contains the shadow of $B(y_\xi,\lambda')$.
Hence Assumption (7) of Theorem \ref{theo:resucborelcantelli}
follows.

\medskip Let $E_f$ be the set of points of $\partial_\infty X_0$
which belong to infinitely many 
$$
B_r(f_3(n_r))=
\N_r\Big(c_3\;e^{-NE[\frac{D(r)}{N}]}\;
f\big(NE[\mbox{$\frac{D(r)}{N}$}]\big)\Big)
$$
for $r$ in $R_0$, and similarly let $E'_f$ be the set of points of
$\partial_\infty X_0$ which belong to infinitely many
$\N_r\big(f(D(r))\;e^{-D(r)} \big)$. As $f$ is slowly varying, there
exists a constant $c_8\geq 1$ such that $E_{\frac{1}{c_8}f}\subset
E'_f\subset E_{c_8f}$. Hence in order to prove Theorem
\ref{theo:maintechnical}, we only have to prove that if
$\int_{1}^{+\infty} f^{\delta-\delta_0}$ converges (resp.~diverges),
then $\mu_{\Ga_0x_0}(E_f)=0$ (resp.~$\mu_{\Ga_0x_0}({}^c E_f)=0$). If
this integral converges, then the result follows from Part [A] of
Theorem \ref{theo:resucborelcantelli}.

\medskip 
In the divergence case, Part [B] of Theorem
\ref{theo:resucborelcantelli} implies that $\mu_{\Ga_0x_0}(E_f)>0$.
Using the ergodicity of $\mu_x$ in a similar way to the end of the
proof of Theorem 5.1 of \cite{HP2}, and the fact that
$\mu_x(\Lambda\Ga_0)=0$ as $\delta_0<\delta$, it follows that $E_f$
has full measure.  
\cqfd

\bigskip Let us now proceed towards our main result, Theorem
\ref{theo:mainmain}.

Let $\epsilon$ be a positive real number and let $g:[0,+\infty[\;\ra
[0,+\infty[$ be a map such that $t\mapsto f(t) = e^{-g(t)}$ is slowly
varying. A geodesic line $\ell$ in $X$ will be called {\it
  $(\epsilon,g)$-Liouville} with respect to $(\Ga,\Ga_0)$ if there
exist a sequence $(t_n)_{n\in \NN}$ of positive times converging to
$+\infty$ and a sequence $(\ga_n)_{n\in \NN}$ of elements of $\Ga$
such that $\ell(t)$ belongs to $\N_\epsilon (\ga_n C_0)$ for every $t$
in $[t_n,t_n+ g(t_n)]$. The following remark implies that up to
changing $g$ by an additive constant (or equivalently up to changing
$f$ by a multiplicative constant), being $(\epsilon,g)$-Liouville does
not depend on $\epsilon$, and depends only on the asymptotic class of
$\ell$.

\brema\label{rem:liouv} 
{\rm (1) Note that if $\epsilon'\geq \epsilon$
  and $g'\leq g$, then a geodesic line which is
  $(\epsilon,g)$-Liouville is $(\epsilon',g')$-Liouville.

(2) Note that by the hyperbolicity properties of $X$, for every
$\epsilon'$ in $]\,0,\epsilon]$, there exists a constant
$c(\epsilon,\epsilon')\geq 0$ such that for every convex subset $C$ of
$X$ and every geodesic line $\ell$ in $X$, if the length $h$ of the
intersection of $\ell$ and $\N_\epsilon C$ is at least
$c(\epsilon,\epsilon')$, then the length of the intersection of $\ell$
and $\N_{\epsilon'} C$ is at least $h-c(\epsilon,\epsilon')$ (see
\cite{PP} for precise estimates). In particular, if $g\geq
c(\epsilon,\epsilon')$, then a geodesic line which is
$(\epsilon,g)$-Liouville is
$(\epsilon',g-c(\epsilon,\epsilon'))$-Liouville.

(3) Recall that two geodesic lines $\ell,\ell'$ in $X$ are {\it
  asymptotic} if $d(\ell(t),\ell')$ (or equivalently
$d(\ell'(t),\ell))$ is bounded (or equivalently tends to $0$) as $t$
tends to $+\infty$. Note that $\ell,\ell'$ are asymptotic if and only
if their points at infinity $\ell(+\infty),\ell'(+\infty)$ are equal.

By the strict convexity of $\epsilon$-neighborhoods of convex subsets
of $X$, if $\ell$ is an $(\epsilon,g)$-Liouville geodesic line, and
$\ell'$ is a geodesic line which is asymptotic to $\ell$, then $\ell'$
is $(2\epsilon,g)$-Liouville, as well as $(\epsilon,g-\eta)$-Liouville
for every constant $\eta>0$ such that $g\geq \eta$.  
}  
\erema

Let $\pi_0:X\cup (\partial_\infty X-\partial_\infty C_0)\ra
X_0\cup\partial_\infty X_0$ be the canonical projection. 
The next lemma shows the relation between the (geometric) 
Liouville property of a geodesic line and the fact that its point 
at infinity belongs to a limsup subset considered in Theorem 
\ref{theo:maintechnical}.

\blemm \label{lem:euivliouborcan} %
There exists $c'''>0$ such that for every geodesic line $\ell$
in $X$ such that $\ell(+\infty)\notin
\bigcup_{\ga\in\Ga}\;\ga\;\partial_\infty C_0$,
\begin{itemize}
\item[(1)] if $\ell$ is $(\epsilon,g)$-Liouville, then the point
  $\pi_0(\ell(+\infty))$ belongs to infinitely many subsets
  $\N_r(c'''\;f(D(r))\;e^{-D(r)})$ for $r$ in $R_0$;
\item[(2)] if $\pi_0(\ell(+\infty))$ belongs to infinitely many
  subsets $\N_r(\frac{1}{c'''}\;f(D(r))\;e^{-D(r)})$ for $r$ in $R_0$,
  then $\ell$ is $(\epsilon,g)$-Liouville.
\end{itemize}
\elemm

\dem 
(1) Assume that $\ell$ is $(\epsilon,g)$-Liouville. Up to
replacing $\epsilon$ by $2\epsilon$ and $\ell$ by an asymptotic line,
as $\ell(+\infty)\notin \partial_\infty C_0$ and by Remark
\ref{rem:liouv} (3), we may assume that $\ell(0)$ is the closest point
in $C_0$ to $\ell(+\infty)$.

\medskip
\noindent
\begin{minipage}{10.1cm}
  ~~~ Let $(t_n)_{n\in \NN}$ be a sequence of positive times,
  converging to $+\infty$. Let $(\ga_n)_{n\in \NN}$ in $\Ga$ be such
  that $\ell(t) \in \N_\epsilon (\ga_n C_0)$ for every $t\in [t_n,t_n+
  g(t_n)]$. As $\Ga$ acts properly on $X$, and as its subgroup $\Ga_0$
  acts cocompactly on $C_0$, the family $(\ga C_0)_{\ga\in
    (\Ga-\Ga_0)/\Ga_0}$ is locally finite. Hence $d(\ga_n C_0,C_0)$
  tends to $+\infty$ as $n\ra+\infty$ (otherwise $\ell(+\infty)$ would
  belong to $\ga\;\partial_\infty C_0$ for some $\ga\in\Ga$).  In
  particular, up to extracting a subsequence, $\ga_n\notin \Ga_0$ and
  with $r_n=[\ga_n]\in R_0$, the $r_n$'s are pairwise distinct.
  Furthermore, we may assume that $\ell$ enters $\N_\epsilon(\ga_n
  C_0)$ at the time $t_n$. Let $[p_n,q_n]$ be the shortest segment
  between $C_0$ and $\ga_n C_0$, with $p_n\in C_0$, so that
  $D(r_n)=d(p_n,q_n)$. Let $x_n$ (resp.~$y_n$) be a point of $\ga_n
  C_0$ such that $d(x_n,\ell(t_n))\leq \epsilon$
  (resp.~$d\big(y_n,\ell(t_n+g(t_n))\big) \leq \epsilon$).
\end{minipage}
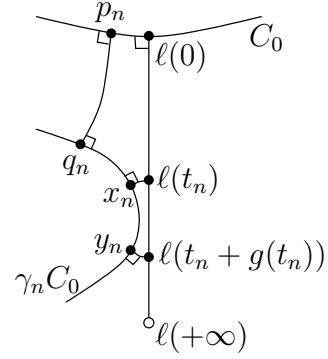
\begin{minipage}{3.8cm}
\begin{center}
\begin{picture}(0,0)%
\includegraphics{fig_equivliou.pstex}%
\end{picture}%
\setlength{\unitlength}{4144sp}%
\begingroup\makeatletter\ifx\SetFigFont\undefined%
\gdef\SetFigFont#1#2#3#4#5{%
  \reset@font\fontsize{#1}{#2pt}%
  \fontfamily{#3}\fontseries{#4}\fontshape{#5}%
  \selectfont}%
\fi\endgroup%
\begin{picture}(1504,2231)(759,-1474)
\put(1613,277){\makebox(0,0)[lb]{\smash{{\SetFigFont{12}{14.4}{\rmdefault}{\mddefault}{\updefault}{\color[rgb]{0,0,0}$\ell(0)$}%
}}}}
\put(2176,393){\makebox(0,0)[lb]{\smash{{\SetFigFont{12}{14.4}{\rmdefault}{\mddefault}{\updefault}{\color[rgb]{0,0,0}$C_0$}%
}}}}
\put(1058,-346){\makebox(0,0)[lb]{\smash{{\SetFigFont{12}{14.4}{\rmdefault}{\mddefault}{\updefault}{\color[rgb]{0,0,0}$q_n$}%
}}}}
\put(1624,-458){\makebox(0,0)[lb]{\smash{{\SetFigFont{12}{14.4}{\rmdefault}{\mddefault}{\updefault}{\color[rgb]{0,0,0}$\ell(t_n)$}%
}}}}
\put(1636,-912){\makebox(0,0)[lb]{\smash{{\SetFigFont{12}{14.4}{\rmdefault}{\mddefault}{\updefault}{\color[rgb]{0,0,0}$\ell(t_n+g(t_n))$}%
}}}}
\put(774,-1055){\makebox(0,0)[lb]{\smash{{\SetFigFont{12}{14.4}{\rmdefault}{\mddefault}{\updefault}{\color[rgb]{0,0,0}$\ga_n C_0$}%
}}}}
\put(1260,574){\makebox(0,0)[lb]{\smash{{\SetFigFont{12}{14.4}{\rmdefault}{\mddefault}{\updefault}{\color[rgb]{0,0,0}$p_n$}%
}}}}
\put(1301,-540){\makebox(0,0)[lb]{\smash{{\SetFigFont{12}{14.4}{\rmdefault}{\mddefault}{\updefault}{\color[rgb]{0,0,0}$x_n$}%
}}}}
\put(1246,-818){\makebox(0,0)[lb]{\smash{{\SetFigFont{12}{14.4}{\rmdefault}{\mddefault}{\updefault}{\color[rgb]{0,0,0}$y_n$}%
}}}}
\put(1616,-1396){\makebox(0,0)[lb]{\smash{{\SetFigFont{12}{14.4}{\rmdefault}{\mddefault}{\updefault}{\color[rgb]{0,0,0}$\ell(+\infty)$}%
}}}}
\end{picture}%

\end{center}
\end{minipage}
\medskip

As we have already seen (in the proof of Lemma
\ref{lem:verifsixoldthree}), there exists a constant $c'''_1>0$ such
that $[q_n,y_n]$ is contained in the $c'''_1$-neighborhood of a
geodesic ray $[q_n,\xi_n[$ with $\xi_n\in\ga_n\partial_\infty C_0$,
and such that for $n$ big enough, $d(p_n,\ell(0))\leq c'''_1$. By
hyperbolicity, the distance between $p_n$ and the closest point of
$C_0$ to $\xi_n$ is at most a constant. By arguments similar to the
ones in the proof of Lemma \ref {lem:verifsixoldthree}, it is easy to
prove that there exists a constant $c'''_2\geq 0$ such that $-\log
d_{C_0}(\xi_n,\ell(+\infty))\geq t_n+g(t_n) -c'''_2$ and (using the
fact that $\ell$ enters in $\N_\epsilon(\ga_n C_0)$ at time $t_n$)
that $|t_n-d(p_n,q_n)|\leq c'''_2$. As $f$ is slowly varying, there
exists a constant $c'''\geq 1$ such that
$d_{C_0}(\xi_n,\ell(+\infty))\leq c'''f(D(r_n))\;e^{-D(r_n)}$. This
proves the first assertion.

\medskip (2) Assume now that there exist a sequence
$(r_n=[\ga_n])_{n\in\NN}$ of pairwise distinct elements in $R_0$ and
$\xi_n\in\ga_n \partial_\infty C_0$ such that
$d_{C_0}(\xi_n,\ell(+\infty)) \leq \frac{1}{c'''}\;f(D(r_n))
\;e^{-D(r_n)}$ for every $n$, for some $c'''\geq 1$ big enough, to be
determined later on. Let us prove that $\ell$ is
$(\epsilon,g)$-Liouville. Up to replacing $\epsilon$ by
$\frac{\epsilon}{2}$ and $\ell$ by an asymptotic line, we may assume
as above that $\ell(0)$ is the closest point in $C_0$ to
$\ell(+\infty)$.

By Lemma \ref{lem:depthsdiscrete}, we have that $D(r_n)= d(C_0,\ga_n
C_0)$ tends to $+\infty$ as $n\ra+\infty$ (hence is positive for $n$
big enough). As above, by hyperbolicity, the closest point of $C_0$ to
$\xi_n$ is at distance at most a constant from the closest point of
$C_0$ to $\ga_n C_0$. By hyperbolicity and the definition of
$d_{C_0}$, there exists a constant $c'''_3\geq 0$ such that between
the times $t=D(r_n)$ and $t=D(r_n)+g(D(r_n))+\log c'''-c'''_3$, the
geodesic ray $\ell$ is at distance at most $c'''_3$ from $\ga_n C_0$.
Hence, as in Remark \ref{rem:liouv} (3), there exists a constant
$c'''_4\geq 0$ such that both $\ell(D(r_n)+c'''_4)$ and
$\ell(D(r_n)+g(D(r_n))+\log c'''-c'''_3-c'''_4)$ are at distance at
most $\epsilon$ from $\ga_n C_0$. Hence if $c'''$ is big enough, by
setting $t_n=D(r_n)+c'''_4$ and as $f$ is slowly varying, the second
assertion follows. \cqfd

\btheo \label{theo:mainmain} %
Let $X$ be a proper ${\rm CAT}(-1)$ geodesic metric space. Let $\Ga$
be a non elementary discrete group of isometries of $X$, with finite
critical exponent $\delta$, of divergent type. Let
$\wt{\mu}_{\mbox{\tiny BM}}$ be its Bowen-Margulis measure. Assume
that ${\rm Card}\;\Ga x \cap B(x,n)\asymp e^{\delta n}$, for some
$x\in X$.  Let $(\Ga_i)_{i\in I}$ be a finite family of almost
malnormal convex-cocompact subgroups of infinite index in $\Ga$ with
critical exponents $(\delta_i)_{i\in I}$.  Let $\delta_0^+=\sup_{i\in
  I}\delta_i$ and $\delta_0^-=\inf_{i\in I}\delta_i$. Let
$g:[0,+\infty[\;\ra \;[0,+\infty[$ be a map such that $t\mapsto
f(t)=e^{-g(t)}$ is slowly varying, and let $\epsilon>0$.

If $\int_{1}^{+\infty} f(t)^{\delta-\delta^+_0}\;dt$ diverges
(resp.~$\int_{1}^{+\infty} f(t)^{\delta-\delta^-_0}\;dt$ converges),
then $\wt{\mu}_{\mbox{\tiny BM}}$-almost every (resp.~no) element of
$\G X$ is $(\epsilon,g)$-Liouville with respect to $(\Ga,\Ga_i)$ for
every (resp.~some) $i\in I$.  
\etheo

\rem (1) The result still holds for a countable family $(\Ga_i)_{i\in
  I}$, under the assumption that $\delta_0^+<\delta$. There are
examples of $\Ga\backslash X$ with $X$ and $\Ga$ as in the above
theorem, such that the upper bound of the critical exponents of the
infinite index subgroups in $\Ga$ is equal to (resp.~is strictly less
than) the critical exponent of $\Ga$, as for instance the closed real
hyperbolic $3$-manifolds fibering over the circle (resp.~the closed
quaternionic hyperbolic manifolds, see for instance \cite{Leu}).

(2) Assume in this remark that $X$ is a Riemannian manifold, that
$\Ga$ is cocompact and torsion free, and that $\Ga_0$ is the
stabilizer of a geodesic line.  This correspond to the hypotheses of
Theorem \ref{theo:introspiralclosgeod} (that appear above it).  Then
there might be a simpler proof using symbolic coding, as indicated to
us by V.~Kleptsyn, using the fact that the geodesic flow of
$\Ga\backslash X$ is then conjugated to a suspension of a Bernoulli
shift. But this requires some serious amount of work, since some
geometric features are difficult, to say the least, to translate by
the coding. In our general situation, no such coding is possible
anyway.

\medskip \dem %
Note that the divergence (resp.~convergence) of the integral in the
statement is unchanged if one replaces $f$ by a scalar multiple of it.
Also recall that $\mu_{x_0}(\bigcup_{i\in I,\;\ga\in \Ga} \ga\Lambda
\Ga_i)=0$, since $\delta_0<\delta$ (see Lemma 
\ref{lem:convcocinfindexpstrinf}).

When the index set $I$ has only one element, the result follows from
Theorem \ref{theo:maintechnical}, by considering the conformal density
$(\mu_x)_{x\in X}$ of dimension $\delta$ for $\Ga$ that is used in the
construction (recalled in Section \ref{sec:backnota}) of
$\wt{\mu}_{\mbox{\tiny BM}}$ (which is ergodic since $\Ga$ is of
divergent type), and by the lemmas \ref{lem:euivliouborcan},
\ref{lem:abscontmes}, \ref{lem:bowmarabscontpatsul}.

Using the fact that finite or countable unions of sets of measure $0$
have measure $0$, the result for general $I$ follows.  
\cqfd

\medskip Using the three examples at the beginning of Section
\ref{sec:spiralgeod}, the theorems \ref{theo:introspiralclosgeod},
\ref{theo:introspiraltotgeod}, \ref{theo:introkleinprecinv} in the
introduction follow.

Similarly, to prove Proposition \ref{prop:introtree} of the
introduction, we apply Theorem \ref{theo:mainmain} to $X=T$ the tree
in the statement of Proposition \ref{prop:introtree}, with $I=\{0\}$
and $\Ga_0$ the stabilizer in $\Ga$ of a geodesic line in $T$ mapping
to the cycle $C$ in the statement of Proposition \ref{prop:introtree}.
We replace the map $g$ in Theorem \ref{theo:mainmain} by $g/L$ for the
map $g$ in Proposition \ref{prop:introtree} (the map $t\mapsto
\exp(-g(t)/L)$ is still slowly varying). When $\Ga$ is cocompact and
torsion-free (anyone of these two assumptions may be non satisfied),
then the symbolic dynamics argument alluded to above works easily, and
gives an alternate proof. But no such coding is easy in general, even
for lattices as simple as ${\rm PSL}_2(\FF_q[X])$ in ${\rm
  PSL}_2\big(\FF_q((X))\big)$, see for instance \cite{BP}.

Many other applications are possible, we will only give the next one.

\medskip 
We refer to \cite{GP} (see also \cite{Bou2,HaP}) for the
definitions and basic properties of an hyperbolic building, which in
particular, when locally finite, is a proper ${\rm CAT}(-1)$ geodesic
metric space. For instance, for every integers $p\geq 5, q\geq 3$, let
$(W_p,S_p)$ be the hyperbolic Coxeter system generated by the
reflections on the sides of a right angled regular real hyperbolic
$p$-gon; Bourdon's building $I_{p,q}$ is (see for instance
\cite{Bou2}) the unique (up to isomorphism) hyperbolic building of
dimension $2$, modeled on $(W_p,S_p)$, and whose links of vertices are
bipartite graphs on $q+q$ vertices. It has a cocompact lattice
$\Ga_{p,q}$ with presentation
$$
\langle s_1,\dots, s_p\;| \;
\forall \; i \in\ZZ/p\ZZ\;\; s_i^q = 1, [s_i,s_{i+1}]=1\rangle\;,
$$
where $s_1, \dots, s_p$ are generators of the pointwise stabilizers of
the $p$ panels of a fundamental chamber $\mathfrak C$ of $I_{p,q}$.
If $q$ is even, let $\Ga_0$ be the subgroup (isomorphic to $W_p$)
generated by the elements $s_i^{\frac{q}{2}}$ for $1\leq i\leq p$,
which is, by the simple transitivity of the action of $\Ga_{p,q}$ on
the set of chambers, the stabilizer of a (unique) apartment
$A_{\mathfrak C}$ in $I_{p,q}$ containing $\mathfrak C$.

\brema \label{rem:haglund}
If $q$ is even, then the subgroup $\Ga_0$ is almost
malnormal in $\Ga_{p,q}$.  
\erema

\dem (F.~Haglund) Let $V$ be the union of the closed chambers of
$I_{pq}$ meeting $A_{\mathfrak C}$, which is invariant by $\Ga_0$.  By
convexity (and arguments as in Poincar\'e's theorem about reflection
groups), the subgroup $H$ of $\Ga$ generated by the pointwise
stabilizers of the edges contained in the boundary of $V$ has $V$ as a
strict fundamental domain, and is normalized by $\Ga_0$. The subgroup
$\Ga'$ of $\Ga$ generated by $H$ and $\Ga_0$, which is isomorphic to
their semi-direct product, has finite index in $\Ga$, since $\Ga$ is
discrete and $\Ga_0$ acts transitively on the chambers of
$A_{\mathfrak C}$.  Let $\Ga''$ be a finite index torsion free
subgroup of $\Ga'$ (which exists for instance since $\Ga_{p,q}$ is
linear, see for example \cite{Kap}).

Let us prove that the stabilizer $S=\Ga''\cap \Ga_0$ of $A_{\mathfrak
  C}$ in $\Ga''$ is malnormal in $\Ga''$, which proves the result.
Assume by absurd that there exists $\ga$ in $\Ga''-S$ and $s$ in
$S-\{e\}$ such that $s$ and $\ga s\ga^{-1}$ preserve $A_{\mathfrak
  C}$. By construction, two distinct translates of $A_{\mathfrak C}$
by elements of $\Ga'$ are disjoint. Hence $\ga^{-1}A_{\mathfrak C}$
and $A_{\mathfrak C}$ are disjoint, and both preserved by $s$. The
(unique) shortest segment between $\ga^{-1}A_{\mathfrak C}$ and
$A_{\mathfrak C}$ is then fixed by $s$, which contradicts the fact
that $\Ga''$ is torsion free.  \cqfd

\bcoro \label{coro:casimmhyp} %
Let $X$ be a locally finite thick hyperbolic building modeled on an
hyperbolic Coxeter system $(W,S)$. Let $\Ga$ be a cocompact lattice in
the automorphism group of $X$ with Bowen-Margulis measure $\mu$.  Let
$A$ be an appartment in $X$ whose stabilizer $\Ga_A$ in $\Ga$ acts
cocompactly on $A$ and is almost malnormal in $\Ga$. Denote by $k\geq
1$ the dimension of $A$ (hence of $X$), and by $\delta$ the Hausdorff
dimension of $\partial_\infty X$ (for any visual distance). Let $f\leq
1$ be a slowly varying map, and $\epsilon>0$.

If $\int_{1}^{+\infty} f(t)^{\delta-k+1}\;dt$ converges
(resp.~diverges), then for $\mu$-almost no (resp.~every) $\ell$ in $\G
X$, there exist positive times $(t_n)_{n\in \NN}$ converging to
$+\infty$ such that $\ell(t)$ belongs to $\Ga \N_\epsilon A$
for every $t$ in $[t_n,t_n-\log
f(t_n)]$.  \ecoro

\dem The apartments in an hyperbolic building are convex (for the
${\rm CAT}(-1)$ metric), hence $\Ga_A$ is convex-cocompact with
critical exponent $k-1$. As $X$ is thick, $\Ga_A$ has infinite index
in $\Ga$. The result follows from Theorem \ref{theo:mainmain} (with
$I$ a singleton).  
\cqfd

\bigskip Let us go back to the general situation of Theorem
\ref{theo:mainmain}. The following result is a logarithm law-type
result for the spiraling of geodesic lines in $\Ga\backslash X$ around
$\Ga_0 \backslash C_0$.  For every $\epsilon>0$ fixed, define the {\it
  penetration map} ${\mathfrak p} ={\mathfrak p}_{\Ga \N_\epsilon C_0}
:\G X\times [0,+\infty[\;\ra[0,+\infty[$ in $\Ga \N_\epsilon C_0$ of
the geodesic lines in $X$, in the following way.  For $(\ell,t)\in\G
X\times [0,+\infty[\,$, if $\ell(t)$ does not belong to $\Ga
\N_\epsilon C_0$, then let ${\mathfrak p}(\ell,t)=0$.  Otherwise, let
${\mathfrak p}(\ell,t)$ be the upper bound of the lengths of the
intervals $I$ in $\RR$ containing $t$ such that there exists $\ga$ in
$\Ga$ with $\ell(I$) contained in $\ga \N_\epsilon C_0$.

\btheo \label{theo:spiralloglaw} %
Let $X$ be a proper CAT$(-1)$ geodesic metric space. Let $\Ga$ be a
non elementary discrete group of isometries of $X$, with finite
critical exponent $\delta$, of divergent type, and let
$\wt{\mu}_{\mbox{\tiny BM}}$ be the Bowen-Margulis measure of $\Ga$.
Assume that ${\rm Card}\;\Ga x \cap B(x,n)\asymp e^{\delta n}$, for
some $x\in X$.  Let $\Ga_0$ be an almost malnormal convex-cocompact
subgroup of infinite index in $\Ga$ with critical exponent $\delta_0$.

Then for every $\epsilon>0$, for $\wt{\mu}_{\mbox{\tiny BM}}$-almost
every $\ell$ in $\G X$, we have
$$
\limsup_{t\ra+\infty}\;
\frac{{\mathfrak p}(\ell,t)}{\log t} 
=\frac{1}{\delta-\delta_0}\;.
$$ 
\etheo

\dem %
For every $\ga$ in $\Ga$ such that a geodesic line $\ell$ enters the
$\epsilon$-neighborhood of $\ga C_0$, let $t_{\ell,\ga}$ be the
entering time of $\ell$ in this neighborhood.

We apply Theorem \ref{theo:mainmain} with $g_\kappa:t\mapsto \kappa
\,\log(1+t)$, which is a Lipschitz map $\RR^+\ra\RR^+$.  Note that the
integral $\int_{1}^{+\infty} t^{-(\delta-\delta_0)\kappa} \;dt$
diverges if and only if $\kappa \leq \frac{1}{\delta-\delta_0}$. If
$\kappa_n =\frac{1}{\delta-\delta_0}+ \frac{1}{n}$ for $n\in
\NN-\{0\}$, then the convergence part of Theorem \ref{theo:mainmain}
implies that for $\wt{\mu}_{\mbox{\tiny BM}}$-almost every $\ell$ in
$\G X$, for every $\ga$ in $\Ga$ such that $\ell$ meets
$\ga\N_\epsilon C_0$ with $t_{\ell,\ga}$ big enough, we have
${\mathfrak p}(\ell,t_{\ell,\ga})\leq g_{\kappa_n}(t_{\ell,\ga})$.
Hence
$$
\limsup_{t\ra+\infty}\; \frac{{\mathfrak p}(\ell,t)}{\log t}=
\limsup\; \frac{{\mathfrak p}(\ell,t_{\ell,\ga})}{\log
  (1+t_{\ell,\ga})}\leq \kappa_n \:,
$$
where the upper limit is taken on the $\ga\in\Ga-\Ga_0$ such that
$\ell$ meets $\ga\N_\epsilon C_0$ and $t_{\ell,\ga}$ tends to
$+\infty$. As $n\ra+\infty$, we get that $\limsup_{t\ra+\infty}\;
\frac{{\mathfrak p}(\ell,t)}{\log t} \leq\frac{1}{\delta-\delta_0}$.
Similarly, using the divergence part of Theorem \ref{theo:mainmain}
with the function $g=g_\kappa$ where $\kappa = \frac{1} {\delta -
 \delta_0}$, we get that for $\wt{\mu}_{\mbox{\tiny BM}}$-almost
every $\ell$ in $\G X$, $\limsup_{t\ra+\infty}\; \frac{{\mathfrak
    p}(\ell,t)}{\log t} \geq\frac{1}{\delta-\delta_0}$.  \cqfd

\medskip 
Corollary \ref{theo:introloglaw} in the introduction follows
immediately.

\section{Non-archimedean Diophantine approximation by 
quadra\-tic irrational numbers}
\label{sec:dioapp}

Let us now give an application of our results to Diophantine
approximation in non-archimedian local fields.

Let $\wh K= \FF_q((X^{-1}))$ be the field of formal Laurent series in
the variable $X^{-1}$ over the finite field $\FF_q$. Recall the
definition of the absolute value of an element $f\in\wh K-\{0\}$. Let
$f=\sum_{i=n}^\infty a_i X^{-i}$ where $n\in \ZZ$ and $a_n\neq 0$.
Then we define $\nu(f)=n$ and $|f|_\infty=q^{-\nu(f)}$. Endow the
locally compact additive group $\wh K$ with its (unique up to a
constant factor) Haar measure $\mu$. Let $K=\FF_q(X)$.

Let $\TT_q$ be the Bruhat-Tits tree of $({{\rm SL}}_2, \wh K)$; we
refer to \cite{Ser} for any background on $\TT_q$.  Identify as usual
$\partial_\infty \TT_q$ and $\wh K\cup\{\infty\}$, so that the action
of ${\rm SL}_2(\wh K)$ on $\TT_q$ extends continuously by the action
by homographies of ${\rm SL}_2(\wh K)$ on $\wh K\cup\{\infty\}$. Let
$x_0$ be the standard base point in $\TT_q$.  Note that the Hausdorff
dimension of the visual distance $d_{x_0}$ is $\log q$, as $\TT_q$ is
a regular tree of degree $q+1$.

We refer for instance to \cite{Las,Sch} for nice surveys of the
Diophantine approximation properties of elements in $\wh K$ by
elements in $K$, a geometric interpretation of which being given in
\cite{Pau}. Here, we are interested in approximating elements of $\wh
K$ by elements in the set $K_2$ of irrational quadratic elements in
$\wh K$ over $K$. For every $\alpha$ in $K_2$, let $\alpha^*$ be its
Galois conjugate (the other root of its minimal polynomial), and
define its {\it height} by
$$
h(\alpha)={|\alpha-\alpha^*|_\infty}^{-1} \;.
$$
We will not make precise here the relationship with the standard 
height of an element of the projective line over the algebraic closure 
of $K$, see for instance \cite{HS}.

Let $\Ga={\rm PSL}_2(\FF_q[X])={\rm SL}_2(\FF_q[X])/\{\pm {\rm id}\}$,
which is a (non-uniform) lattice of $\TT_q$ (see for instance
\cite{Ser}), hence a non-elementary discrete group of isometries of
$\TT_q$, whose critical exponent $\delta$ is equal to the Hausdorff
dimension of $d_{x_0}$, that is $\delta=\log q$. See for instance
\cite{BP} for a (well known) proof that the restrictions to
$\partial_\infty \TT_q-\{\infty\}=\wh K$ of the Patterson-Sullivan
measures of $\Ga$ have the same measure class as the Haar measure
$\mu$ of $\wt K$.

\medskip
\noindent{\bf Proof of Theorem \ref{theo:intrononarchapproxdioph}. }
Let $\ga_0$ be an hyperbolic element of $\Ga$, $C_0$ be its
translation axis in $\TT_q$, and $\Ga_0$ be the stabilizer of $C_0$ in
$\Ga$, which is convex-cocompact with critical exponent $\delta_0=0$.
It is easy to verify that the set of points at infinity of $C_0$ is
$\{\alpha,\alpha^*\}$ for some $\alpha$ in $K_2$; and that any such
pair is the set of endpoints of some hyperbolic element of $\Ga$ (one
can for instance use the fact that Artin's continued fraction
expansion of an element in $K_2$ is eventually periodic (see for
example \cite{Las})).

Note that for every $\ga\in\Ga$ and $\alpha\in K_2$, the element $\ga
\alpha$ is still in $K_2$, $(\ga\alpha)^*=\ga\alpha^*$ and
$\ga\{\alpha,\alpha^*\}\cap\{\alpha,\alpha^*\}\neq\emptyset$ if and
only if $\ga\in\Ga_0$.

Denote by $d_\infty$ the Hamenst\"adt distance on
$\partial_\infty\TT_q- \{\infty\} =\wh K$ defined by the horosphere
centered at $\infty$ and passing through $x_0$. It is proved in
\cite[Coro.~5.2]{Pau} that
$d_\infty(\xi,\xi')=|\xi-\xi'|_\infty^{\frac{1}{\log q}}$, for every
$\xi,\xi'$ in $\wh K$.

\blemm For every $\xi_0$ in $\wh K-\{\alpha,\alpha_*\}$, there exists
a neighborhood $V$ of $\xi_0$ and a constant $c_*>0$ such that for
every $\xi,\xi'$ in $V$,
$$
d_{C_0}(\xi,\xi')=c_*\;d_\infty(\xi,\xi')=
c_*\;|\xi-\xi'|_\infty^{\frac{1}{\log q}}\;.
$$
Furthermore, for every $\ga$ in $\Ga-\Ga_0$ such that $\ga\alpha$ and
$\ga\alpha^*$ belongs to $V$, we have
$$
e^{-D([\ga])}=c_*\;|\ga\alpha-\ga\alpha^*|_\infty^{\frac{1}{\log q}}\;.
$$
\elemm

\medskip
\noindent
\begin{minipage}{10.6cm}
  \dem Let $p_0$ be the intersection of the geodesic line
  $]\infty,\xi_0[$ in $\TT_q$ with the horosphere centered at $\infty$
  passing through $x_0$. Let $q_0=\pi_{C_0}(\xi_0)$, and $u_0\in X$
  such that $]\xi_0,\infty[\;\cap\; ]\xi_0,q_0]= \;]\xi_0,u_0]$.  In
  the picture on the right, we assume that $C_0$ and $]\infty,\xi_0[$
  are disjoint, and that $p_0\in[u_0,\infty[$. But the following
  reasoning is independent of these assumptions. Let $c_*=
  e^{-\beta_{\xi_0}(q_0,p_0)}$. If $\xi,\xi'$ are close enough to
  $\xi_0$, and $]\xi,\infty[\;\cap\;]\xi',\infty[\;= [p,\infty[\,$,
  then $p_0,u_0\in[p,\infty[\,$, $\pi_{C_0}(\xi)=\pi_{C_0}(\xi')=q_0$,
  $\beta_{\xi_0}(q_0,p_0)=d(q_0,p)-d(p_0,p)$, $d_\infty(\xi,\xi')=
  e^{-d(p_0,p)}$ and $d_{C_0}(\xi,\xi')= e^{-d(q_0,p)}$, hence the
  first result follows.

  ~~~ As $d_{C_0}(\ga\alpha,\ga\alpha^*)=e^{-D([\ga])}$ if $\ga
  \alpha,\ga\alpha^*$ are closed enough to $\xi_0$ (see Equation
  \eqref{eq:calcdcarb}), the second result follows from the first one.
  \cqfd
\end{minipage}
\begin{minipage}{4.2cm}
\begin{center}
\begin{picture}(0,0)%
\includegraphics{fig_techarb.pstex}%
\end{picture}%
\setlength{\unitlength}{3398sp}%
\begingroup\makeatletter\ifx\SetFigFont\undefined%
\gdef\SetFigFont#1#2#3#4#5{%
  \reset@font\fontsize{#1}{#2pt}%
  \fontfamily{#3}\fontseries{#4}\fontshape{#5}%
  \selectfont}%
\fi\endgroup%
\begin{picture}(1932,2466)(1321,-2749)
\put(2161,-466){\makebox(0,0)[lb]{\smash{{\SetFigFont{10}{12.0}{\rmdefault}{\mddefault}{\updefault}{\color[rgb]{0,0,0}$\infty$}%
}}}}
\put(2611,-2581){\makebox(0,0)[lb]{\smash{{\SetFigFont{10}{12.0}{\rmdefault}{\mddefault}{\updefault}{\color[rgb]{0,0,0}$\xi$}%
}}}}
\put(2071,-2671){\makebox(0,0)[lb]{\smash{{\SetFigFont{10}{12.0}{\rmdefault}{\mddefault}{\updefault}{\color[rgb]{0,0,0}$\xi'$}%
}}}}
\put(2881,-1546){\makebox(0,0)[lb]{\smash{{\SetFigFont{10}{12.0}{\rmdefault}{\mddefault}{\updefault}{\color[rgb]{0,0,0}$q_0$}%
}}}}
\put(2259,-2131){\makebox(0,0)[lb]{\smash{{\SetFigFont{10}{12.0}{\rmdefault}{\mddefault}{\updefault}{\color[rgb]{0,0,0}$p$}%
}}}}
\put(3098,-2041){\makebox(0,0)[lb]{\smash{{\SetFigFont{10}{12.0}{\rmdefault}{\mddefault}{\updefault}{\color[rgb]{0,0,0}$C_0$}%
}}}}
\put(1336,-2671){\makebox(0,0)[lb]{\smash{{\SetFigFont{10}{12.0}{\rmdefault}{\mddefault}{\updefault}{\color[rgb]{0,0,0}$\xi_0$}%
}}}}
\put(1876,-1659){\makebox(0,0)[lb]{\smash{{\SetFigFont{10}{12.0}{\rmdefault}{\mddefault}{\updefault}{\color[rgb]{0,0,0}$u_0$}%
}}}}
\put(2535,-879){\makebox(0,0)[lb]{\smash{{\SetFigFont{10}{12.0}{\rmdefault}{\mddefault}{\updefault}{\color[rgb]{0,0,0}$x_0$}%
}}}}
\put(2236,-1246){\makebox(0,0)[lb]{\smash{{\SetFigFont{10}{12.0}{\rmdefault}{\mddefault}{\updefault}{\color[rgb]{0,0,0}$p_0$}%
}}}}
\end{picture}%

\end{center}
\end{minipage}
\medskip

\medskip Let $\varphi:[0,+\infty[\;\ra\;]0,1]\;$ be a map with
$t\mapsto f(t)=\varphi (q^t)^{\frac{1}{\log q}}$ slowly varying, and
let $g:t\mapsto -\log f(t)= -\log_q\varphi(q^t)$, so that $\varphi(t)=
q^{- g(\log_q t)}$. By an easy change of variable, the
integral $\int_{1}^{+\infty} \varphi(t)/t \;dt$ diverges if and only
if $\int_{1}^{+\infty} f(t)^{\log q}\,dt$ diverges.

By the above lemma and as $f$ is slowly varying, for every compact
subset $A$ of $\wh K-\{\alpha,\alpha_*\}$, there exist positive
constants $c'_*,c''_*$ such that for every $\xi$ in $A$, 
\begin{itemize}
\item if $(r_n=[\ga_n])_{n\in\NN}$ is a sequence in $R_0$ with
  $D(r_n)\ra+\infty$ as $n\ra+\infty$ and
  $d_{C_0}(\xi,\ga_n\alpha)$ $\leq c'_*\;f(D(r_n))\;e^{-D(r_n)}$ for
  every $n$ big enough, then $h(\ga_n\alpha)\ra+\infty$ as
  $n\ra+\infty$ and, for every $n$ big enough,
$$
|\xi-\ga_n\alpha|_\infty^{\frac{1}{\log q}}\leq
\; e^{-g(-\log |\ga_n\alpha-\ga_n\alpha^*|_\infty^{\frac{1}{\log q}})}\;
  |\ga_n\alpha-\ga_n\alpha^*|^{\frac{1}{\log q}}\;,
$$
that is
$$
|\xi-\ga_n\alpha|_\infty\leq
\frac{\varphi(h(\ga_n\alpha))}{h(\ga_n\alpha)}\;;
$$
\item conversely, if $(\ga_n)_{n\in\NN}$ is a sequence in $\Ga$ with
  $h(\ga_n\alpha)\ra+\infty$ as $n\ra+\infty$ (in particular
  $\ga_n\notin \Ga_0$ for $n$ big enough) and $|\xi- \ga_n\alpha|
  _\infty \leq \frac{\varphi(h(\ga_n\alpha))} {h(\ga_n\alpha)}$ for
  every $n$ big enough, then with $r_n=[\ga_n]$, we have
  $D(r_n)\ra+\infty$ as $n\ra+\infty$ and 
  $d_{C_0}(\xi,\ga_n\alpha) \leq c''_*\; f(D(r_n))\;e^{-D(r_n)}$
  for every $n$ big enough.
\end{itemize}

Hence by Theorem \ref{theo:maintechnical}, if $\int_{1}^{+\infty}
\varphi(t)/t \;dt$ diverges, then for $\mu$-almost every $\xi$ in $\wh
K$, there exist a sequence $(\beta_n)_{n\in\NN}$ in the congruence
class of $\alpha$ in $K_2$, with $h(\beta_n)\ra+\infty$ as
$n\ra+\infty$, such that $|\xi-\beta_n|_\infty\leq
\frac{\varphi(h(\beta_n))}{h(\beta_n)}$.  This can be written as
$\liminf \;\frac{\varphi(h(\beta))}{h(\beta)}\;|\xi-\beta|_\infty\leq
1$, where the lower limit is taken over $\beta$ in the congruence
class of $\alpha$ with $h(\beta)\ra +\infty$. Replacing $\varphi$ by
$\frac{1}{k}\varphi$ and letting $k$ go to $+\infty$, this proves the
divergence part of Theorem \ref{theo:intrononarchapproxdioph} in the
introduction. The convergence part follows similarly. \cqfd

\medskip 
By taking $\varphi:t\mapsto t^{-s}$ in Theorem
\ref{theo:intrononarchapproxdioph} with $s\geq 0$, the next result,
which in particular says that almost every element of $\wh K$ is badly
approximable by quadratic irrational elements of $\wh K$, follows
immediatly.

\bcoro 
For $\mu$-almost every $x$ in $\wh K$, $\liminf h(\beta)
|x-\beta|_\infty=0$, and, for every $s>0$, $\lim
h(\beta)^{1+s}|x-\beta|_\infty=+\infty$, where the lower limit and
limit are taken over the quadratic irrational elements $\beta$ in $\wh
K$, in any (resp.~some) congruence class, with $h(\beta)\ra +\infty$.
\cqfd \ecoro

\section{Approximating points} 
\label{sec:approxpoint}

Let $X$ be a proper ${\rm CAT}(-1)$ geodesic metric space. Let $\Ga$
be a non elementary discrete group of isometries of $X$, with finite
critical exponent $\delta$. In this section (which is a joint work
with C.~S.~Aravinda), we will also apply our
geometric avatar of the Borel-Cantelli Lemma, Theorem
\ref{theo:resucborelcantelli}, to prove a Khintchine-type result for
the approximation of a point by geodesic lines in $X$.

Let $x_0\in X$ be a base point. For every $C\geq 0$, a point $\xi$ in
$\partial_\infty X$ will be called a {\it $C$-strongly conical limit
point} if there exist a geodesic line $\rho$ with $\rho(+\infty)=\xi$
and a sequence $(\ga_n)_{n\in\NN}$ in $\Ga$, such that 
$(\ga_n x_0)_{n\in\NN}$ converges to $\xi$, $d(\ga_n
x_0,\rho)\leq C$ and $d(\ga_n x_0,\ga_{n+1}x_0)\leq C$.  Note that if
$\xi$ is a $C$-strongly conical limit point with respect to $x_0$,
then $\xi$ is a $C'$-strongly conical limit point with respect to any
other base point $x'_0$ for $C'=C+2\,d(x_0,x'_0)$. And if $\xi$ is a
$C$-strongly conical limit point for the geodesic line $\rho$, then
$\xi$ is a $(C+\epsilon)$-strongly conical limit point with respect to
any other geodesic line $\rho'$ asymptotic to $\rho$, for every
$\epsilon>0$.

\medskip \noindent {\bf Examples. } 

(1) If $\xi$ is a fixed point of an hyperbolic element $\ga$ of $\Ga$,
then $\xi$ is a $C$-strongly conical limit point with
$C=\max\{d(x_0,A_\ga), d(x_0,\ga x_0)\}$, where $A_\ga$ is the
translation axis of $\ga$.

(2) If $\Ga$ is convex-cocompact, then there exists a constant $C\geq
0$ such that any limit point of $\Ga$ is a $C$-strongly conical limit
point.

\medskip The following result is (a slight adaptation of) the
well-known Sullivan's shadow lemma, see for instance \cite[page
93]{Bou}.

\blemm \label{lem:sullivanshadow} For every conformal density
$(\mu_z)_{z\in X}$ of dimension $\delta$ for $\Ga$, for every $C\geq
0$, there exists $c\geq 1$ such that for every $C$-strongly conical
limit point $\xi$, for every $\epsilon\in\;]\,0,1]$,
$$
\frac{1}{c} \;\epsilon^\delta\leq 
\mu_{x_0}(B_{d_{x_0}}(\xi,\epsilon))\leq 
c\;\epsilon^\delta\;.\;\;\;\mbox{\cqfd}
$$
\elemm

The following result is the main technical tool from which
Theorem \ref{theo:introapproxpoint} in the introduction follows.

\btheo\label{theo:mainapproxpointtechnical} Let $X$ be a ${\rm
  CAT}(-1)$ proper geodesic metric space. Let $x,y$ be points in $X$.
Let $\Ga$ be a non elementary discrete subgroup of isometries of $X$.
Let $(\mu_z)_{z\in X}$ be a
conformal density of dimension $\delta$ for $\Ga$, for some $\delta$
in $]0,+\infty[\,$. Assume that there exist constants $c_0,C_0>0$ such
that ${\rm ~Card~}\{\ga\in\Ga \;:\; d(x,\ga y)\leq n\}\sim
c_0\;e^{\delta n}$, and that for every $z$ in $\Ga y$ except finitely
many of them, there exists a geodesic ray $\rho_{z}$ starting from
$x$, passing through $z$ and ending at a $C_0$-strongly conical limit
point. Let $f:[\,0,+\infty[\;\ra \;]\,0,+\infty\,[\,$ be a slowly
varying map, with $f(t)$ converging to $0$ as $t\ra+\infty$.  Let
$E_f$ be the set of points in $\partial_\infty X$ which belong to
infinitely many balls $B_{d_x}\big(\rho_{\ga y}(+\infty), f(d(x,\ga
y))\;e^{-d(x,\ga y)}\big)$ for $\ga$ in $\Ga$.

\smallskip 
[A] If $\int_{1}^\infty f(t)^\delta dt$ converges, then
$\mu_x(E_f)=0$.

\smallskip [B] If there exists a sequence $(t_n)_{n\in\NN}$ in $\RR$,
with $t_n\ra+\infty$ as $n\ra+\infty$ and $t_{n+1}\geq t_n -\log
f(t_n)$ for every $n$, such that $\sum_{n\in\NN} f(t_n)^\delta$
diverges, then $\mu_x(E_f)>0$.  
\etheo

\dem 
Since $\Ga$ is discrete, the stabilizer $\Ga_y$ of $y$ in $\Ga$
is finite, and there exists $c_1\in\;]\,0,1]$ such that for every
$\ga,\ga'$ in $\Ga$, if $d(\ga y,\ga' y)\leq 4\,c_1$, then $\ga y=\ga'
y$, that is the class of $\ga$ and $\ga'$ in $\Ga/\Ga_y$ are equal.
It can be easily checked that for $[\ga]$ in $\Ga/\Ga_y$, the objets
$\rho_{\ga y}$ and $d(x,\ga y)$ are well defined.

Let $t'_0> 0$ be big enough so that $f(t)\leq \frac{1}{e}$ for $t\geq
t'_0$, and that for every $\ga$ in $\Ga$ such that $d(x,\ga y) \geq
t'_0$, the geodesic ray $\rho_{\ga y}$ is defined. In case [A], define
$t_0=t'_0$ and by induction $t_{n+1}=t_n+1$ for every $n$ in $\NN$. In
case $[B]$, as $t_n\ra+\infty$, we may assume, up to shifting the
indices, that $t_0\geq t'_0$.

For every $n$ in $\NN$, let
$$
I_n=\{\ga\in \Ga/\Ga_{y}\;:\; t_n\leq d(x,\ga y)<
t_n +c_1\}\;,
$$ 
and $I=\bigcup_{n=0}^{\infty}I_n$. By the discreteness of $\Ga$, the
subsets $I_n$ are finite, and pairwise disjoint since $c_1\leq 1$,
using also the hypothesis on the range of $f$ in case [B].

Since $f$ is slowly varying, there exists $c_2\geq 1$ such that
$f(y)\leq c_2f(x)$ if $|y-x|\leq c_1$. Let $g=-\log f$, which is non
negative on $[t'_0,+\infty[\,$, and $c_3=\min\{c_1,\frac{1}{c_2}\}
\,e^{-c_1}$, whicqh belongs to $]\,0,1]$.  Note that the constants
$c_1,c_2$, hence $c_3$, are unchanged if one replaces $f$ by a scalar
multiple of it.  For every $n$ in $\NN$ and $\epsilon>0$, define
$$
f_1(n)=e^{\delta t_n},\;\;
f_2(n)=c_3\,e^{- t_n},\;\;
f_3(n)=c_3\,e^{-t_n -g(t_n)},\;\;
f_4(n)=c_3^{-\delta},\;\;f_5(\epsilon)=\epsilon^{\delta}\;.
$$ 
The series $\sum_{n\in\NN} f_1(n)f_4(n)f_5(f_3(n))$
converges if and only if the series $\sum_{n\in\NN} f(t_n)^\delta$
converges, which in case [A] is true if and only if the integral
$\int_{1}^{+\infty} f^{\delta}$ converges, as $f$ is slowly varying.

For every $\ga$ in $I$ and $\epsilon>0$, let $B_\ga(\epsilon)=
B_{d_x}(\rho_{\ga y}(+\infty), \epsilon)$, which is measurable and
non-decreasing in $\epsilon$.

\medskip 
Let us prove that the finite measured space $(\partial_\infty
X, \mu_x)$, the family $(B_\ga(\epsilon))_{\ga\in I,\,\epsilon>0}$,
the finite-to-one map $I\ra\NN$ defined by $\ga\mapsto n_\ga=n$ if
$\ga\in I_{n}$, and the above maps $f_1, f_2, f_3, f_4$, $f_5$ satisfy
the assumptions of Theorem \ref{theo:resucborelcantelli}.

Assumption (1) of Theorem \ref{theo:resucborelcantelli} is satisfied
as $g(t_n)\geq 0$ for every $n$. Assumption (2) is satisfied by the
definition of $f_4$. Assumption (3) holds true with for instance
$c'=2,c''=2^\delta$. Assumption (4) follows by the hypothesis on the
growth of the orbit $\Ga y$ (counting elements of $\Ga/\Ga_y$ rather
than of $\Ga$ only divides the number by ${\rm Card}(\Ga_y)$ ).
Assumption (5) is satisfied by Lemma \ref{lem:sullivanshadow}, as
$f_2\leq 1$ and $f_4$ is constant.

\medskip Let us prove Assumption (6) of Theorem
\ref{theo:resucborelcantelli}. Given $n$ in $\NN$ and distinct
$\alpha, \beta$ in $I_n$, assume by absurd that $B_\alpha(f_2(n))$
intersects $B_\beta(f_2(n))$ non trivially. By Lemma
\ref{lem:HPcrelle} (1), the ball $B_{d_x}(\rho_{\ga y}(+\infty),
c_1\;e^{-d(x,\ga y)})$ is contained in $\O_x B(\ga\, y,c_1)$ for every
$\ga\in I$. As $d(x,\alpha\, y) \leq t_n+c_1$ and $c_3\leq c_1\,e^{-
  c_1}$, we have $f_2(n)\leq c_1\,e^{-d(x,\alpha\, y)}$, and similarly
for $\beta$. Hence $\O_x B(\alpha\, y,c_1)$ intersects $\O_x B(\beta\,
y,c_1)$ non trivially.  There exists a geodesic ray $\rho$ starting
from $x$ and passing at distance at most $c_1$ from both $\alpha\, y$
and $\beta\, y$. Let $p,q$ be the closest point of $\alpha y,\beta y$
respectively on $\rho$, with (up to permuting $\alpha$ and $\beta$)
$q\in [x,p]$. As $\alpha,\beta\in I_n$ and since the closest point
maps do not increase distances, we have
$$
d(\alpha\, y,\beta \,y)\leq d(p,q)+2\,c_1 =d(x,p)-d(x,q)+2\,c_1\leq 
d(x,\alpha y)-d(x,\beta_y)+3\,c_1 \leq 4\,c_1\;.
$$ 
This contradicts the definition of $c_1$.

\medskip Finally, let us prove Assumption (7) of Theorem
\ref{theo:resucborelcantelli} under the hypotheses of Case [B].  For
$m>n$, take $(\alpha,\beta)$ in $I_n\times I_m$, and assume that
$B_\beta(f_3(m))$ intersects $B_\alpha(f_3(n))$ non trivially. Let $w$
be a common point of these balls. Since $f_3\leq f_2$, $t_m\geq
t_{n+1}$ and $t_{n+1}\geq t_n+g(t_n)$, we have $f_3(m)\leq f_2(m)\leq
f_2(n+1)\leq f_3(n)$.  Then, for every $z$ in $B_\beta(f_2(m))$, we
have
\begin{align*}
d(z,\rho_{\alpha y}(+\infty))&\leq
d(z,\rho_{\beta y}(+\infty))+d(\rho_{\beta y}(+\infty),w)+
d(w,\rho_{\alpha y}(+\infty))\\ & \leq f_2(m) + f_3(m)+ f_3(n)
\leq 3f_3(n)\;.\end{align*}
Therefore $B_\beta(f_2(m))$ is contained in $B_\alpha(3f_3(n))$, which
proves the claim.

\medskip Let $E'_f$ be the subset of points of $\partial_\infty X$
which, as $\ga$ ranges over $I$, belong to infinitely many balls
$B_\ga(f_3(n_\ga))$. As $0\leq d(x,\ga y)-t_n\leq c_1$ if $\ga\in
I_n$, and since $c_3\leq \frac{1}{c_2}\,e^{-c_1}$, we have, for every
$\ga$ in $\Ga$,
\begin{equation}\label{eq:bouboule}
B_\ga(f_3(n_\ga))\subset
B_{d_x}\big(\rho_{\ga y}(+\infty), f(d(x,\ga y))\;e^{-d(x,\ga y)}\big)
\subset B_\ga\Big(\frac{c_2}{c_3}f_3(n_\ga)\Big)
\;.
\end{equation}

In case [B], it follows from Theorem \ref{theo:resucborelcantelli}
that $\mu_x(E'_f)>0$. By Equation \eqref{eq:bouboule}, we have
$\mu_x(E_f)>0$.

In case [A], since the convergence or divergence of the integral
$\int_1^{+\infty} f^\delta$ is unchanged if one replaces $f$ by a
scalar multiple of it, Theorem \ref{theo:resucborelcantelli} implies
that $\mu_x(E'_{\frac{c_2}{c_3}f})=0$. It follows similarly from
Equation \eqref{eq:bouboule} that $\mu_x(E_f)=0$.  
\cqfd

\medskip Given a complete Riemannian manifold $M$, recall (see
\cite{HP1}) that an element $v$ in $T^1M$, or its associated geodesic
line $\rho$ in $M$, is called {\it $f$-Liouville} at $x_0$ if there
exists a sequence of times $(t_n)_{n\in\NN}$ tending to $+\infty$ such
that $d(\rho(t_n),x_0)\leq f(t_n)$ for every $n$.
The following result is joint work with C.~S.~Aravinda.

\btheo\label{theo:introapproxpoint} Let $M$ be a closed manifold with
sectional curvature at most $-1$, let $x_0$ be a point in $M$, and let
$\mu$ be the maximal entropy probability measure for the geodesic flow
of $M$, with $h$ its topological entropy. Let $f$ be a slowly varying
map.

(1) If $f(t)$ and $e^{-t}/f(t)$ converge to $0$ as $t\ra+\infty$, 
if the integral $\int_u^\infty \frac{f(t)^{h}}{-\log f(t)} \;dt$
diverges (for some $u$ big enough), then $\mu$-almost every geodesic
line is $f$-Liouville at $x_0$.

(2) If the sectional curvature of $M$ satisfies $-a^2\leq K\leq -1$,
and if $\int_1^\infty f(t)^{\frac{h}{a}}\;dt$ converges, then
$\mu$-almost no geodesic line is $f$-Liouville at $x_0$.  \etheo

\dem Let $X\ra M$ be a universal covering of $M$, with covering group
$\Ga$ (which is non elementary), let $x$ be a lift of $x_0$, and take
$y=x$. For every $z$ in $\Ga y-\{y\}$, let $\rho$ be the geodesic ray
starting from $x$ through $z$, which ends at a (uniformly) strongly
conical limit point, as $M$ is compact. Let $(\mu_z)_{z\in X}$ be the
(unique up to positive scalar multiple) ergodic conformal density of
dimension equal to the critical exponent $\delta$ of $\Ga$ (which is
equal to $h$, as $M$ is compact). Since $M$ is compact, we have ${\rm
  ~Card~}\{\ga\in\Ga \;:\; d(x,\ga y)\leq n\}\sim c_0\;e^{\delta n}$
for some constant $c_0>0$. Hence the general hypotheses of Theorem
\ref{theo:mainapproxpointtechnical} on $X,\Ga,(\mu_z)_{z\in X}$ are
satisfied.

\medskip Let us prove the first assertion of Theorem
\ref{theo:introapproxpoint}. Let $f_*:[0,+\infty[\ra \;
]0,\frac{1}{e}]$ be a slowly varying non increasing map such that
$\frac{f_*(t)}{f(t)}$ tends to $0$ as $t\ra+\infty$ (so that in
particular $f_*$ converges to $0$ at $+\infty$) and $\int_u^\infty
\frac{f_*(t)^{h}}{-\log f_*(t)} \;dt$ still diverges.  Fix $t_0\geq 1$
such that $f_*(t_0)\leq \frac{1}{e}$ for $t\geq t_0$.  Define by
induction $t_{n+1}=t_n-\log f_*(t_n)$. In particular, the sequence
$(t_n)_{n\in\NN}$ converges to $+\infty$. As $f_*$ is non increasing,
we have
$$\int_{t_n}^{t_{n+1}} \frac{f_*(t)^\delta}{-\log f_*(t)}\,dt \leq
\frac{f_*(t_n)^\delta}{-\log f_*(t_n)}\,(t_{n+1}-t_n)=f_*(t_n)^\delta\;. 
$$
Hence $\sum_{n\in\NN}f_*(t_n)^\delta$ diverges. Therefore the
hypotheses of Theorem \ref{theo:mainapproxpointtechnical} [B] are
satisfied for $f_*$.

\smallskip Denote by $\S$ the measurable set of elements $\xi$ in
$\partial_\infty X$ which belong to infinitely many visual balls
$B_{d_x}\big(\rho_{\ga x}(+\infty), f_*(d(x,\ga x))\;e^{-d(x,\ga
  x)}\big)$ as $\ga$ ranges over $\Ga$. By Lemma \ref{lem:HPcrelle}
(1), this ball is contained in the shadow $\O_x\big(B(\ga
x,f_*(d(x,\ga x)))\big)$. For every $\xi$ in $\S$, let $\rho_\xi$ be
the geodesic ray starting from $x$ and ending at $\xi$.  As
$f_*(t)\leq 1$ if $t\geq t_0$ and as $f_*$ is slowly varying, there
exist a constant $c_1>0$ and sequences $(s_n)_{n\in\NN}$ in 
$[0,+\infty[$ converging to
$+\infty$ and $(\ga_n)_{n\in\NN}$ in $\Ga$ such that
$d(\rho_\xi(s_n),\ga_n x)\leq c_1\,f_*(s_n)$ for every $n$. For every $\xi$
in $\S$, let $\S_\xi$ be the set of elements $v$ in $T^1M$ such that
the point at infinity of some geodesic line $\rho_v$, which is a lift
by $X\ra M$ of the geodesic line in $M$ associated to $v$, is equal to
$\xi$. Since $X$ is ${\rm CAT}(-1)$, if $\rho,\rho'$ are
asymptotic geodesic rays, then there exists $c>0$ and $\tau\in\RR$
such that $d(\rho(t),\rho'(t+\tau))\leq c\,e^{-t}$ for every $t\geq
\max\{0,-\tau\}$.  Since $f$ is slowly varying, for every $v\in
\S_\xi$, there exist hence a constant $c'>0$ and sequences
$(s'_n)_{n\in\NN}$ in $[0,+\infty[$ converging to $+\infty$ and
$(\ga_n)_{n\in\NN}$ such that
$$
d(\rho_v(s'_n),\ga_n x)\leq c'(f_*(s'_n)+e^{-s'_n})\leq f(s'_n)
\;,
$$
for $n$ big enough. Hence every $v\in \S_\xi$ is $f$-Liouville, for
every $\xi \in \S$. By Theorem \ref{theo:mainapproxpointtechnical}
[B], the set $\S$ has positive measure for $\mu_x$. Hence the set
$\S'=\bigcup_{\xi\in\S}\S_\xi$ (which is measurable) has positive
measure for $\mu$, by Lemma \ref{lem:bowmarabscontpatsul}. As $\S'$ is
invariant under the geodesic flow by construction, and by ergodicity,
it has full measure. This proves the result.

\medskip Let us now prove the assertion (2) of Theorem
\ref{theo:introapproxpoint}. As $f$ is slowly varying, the convergence
of the integral $\int_1^\infty f^{\frac{\delta}{a}}$ implies that $f$
converges to $0$ as $t\ra +\infty$. For every $v$ in $T^1M$, let
$\rho_v$ be a lift by $X\ra M$ of the geodesic line in $M$ associated
to $v$, let $\xi_v=\rho_v(+\infty)$ and let $\rho'_v$ be the geodesic
ray from $x$ to $\xi_v$. If $v$ is $f$-Liouville at $x_0$, then there
exist sequences $(s_n)_{n\in\NN}$ in $[0,+\infty[$ converging to
$+\infty$ and $(\ga_n)_{n\in\NN}$ in $\Ga$ such that
$d(\rho_\xi(s_n),\ga_n x)\leq f(s_n)$ for every $n$. Since $f$ is
slowly varying and since $\rho'_v$ and $\rho_v$ are asymptotic, as
above, there exist $k$ in $\NN-\{0\}$, a sequence $(s'_n)_{n\in\NN}$
in $[0,+\infty[$ converging to $+\infty$, and $(\ga_n)_{n\in\NN}$ in
$\Ga$, such that $d(\rho'_v(s'_n),\ga_n x)\leq k(f(s'_n)+e^{-s'_n})$
for every $n$. In particular, as $f(t)\leq 1$ if $t$ is big enough,
if $v$ is $f$-Liouville, then there exists $k'$ in $\NN-\{0\}$ such
that $\xi_v$ belongs to infinitely many shadows $\O_x\Big(B\big(\ga
x,k'(f(d(x,\ga x))+e^{-d(x,\ga x)})\big)\Big)$ for $\ga\in\Ga$. By
Lemma \ref{lem:HPcrelle} (2), this shadow is contained, except for
finitely many $\ga\in\Ga$, in the ball
$$
\B_{\ga,k''}= B_{d_x}\Big(\rho_{\ga x}(+\infty), k''\big(f(d(x,\ga x))+ 
e^{-d(x,\ga x)}\big)^{\frac{1}{a}}\;e^{-d(x,\ga x)}\Big)\;,
$$
for some positive integer $k''$. If $u,v,w>0$, recall that
$(u+v)^w\leq 2^w(u^w+v^w)$. Hence since $\int_1^\infty
f^{\frac{\delta}{a}}$ converges, the integral $\int_1^\infty \big(
k''(f(t)+ e^{-t})^{\frac{1}{a}}\big)^{\delta}\;dt$ also converges. The
map $t\mapsto k''(f(t)+e^{-t})$ is slowly varying. By Theorem
\ref{theo:mainapproxpointtechnical} [A], the measure of the set of
points in $\partial_\infty X$ which belong to infinitely many balls
$\B_{\ga,k''}$, as $\ga$ ranges over $\Ga$, has measure $0$ for
$\mu_x$. By Lemma \ref{lem:bowmarabscontpatsul}, and since a countable
union of measure zero subsets is a measure zero subset, the result
follows.  
\cqfd

\medskip For every $\alpha>0$, let $f_\alpha:t\mapsto
\frac{1}{(2+t)^{\alpha}}$, which is slowly varying, with $f(t)$ and
$e^{-t}/f(t)$ converging to $0$ as $t\ra+\infty$.  For every $h>0$,
the integral $\int_1^{\infty} \frac{f_\alpha^h}{-\log f_\alpha}$
diverges if and only if $\alpha\leq \frac{1}{h}$ and the integral
$\int_1^{\infty} f_\alpha^h$ converges if and only if
$\alpha>\frac{1}{h}$. By applying Theorem \ref{theo:introapproxpoint}
with $M$ having constant curvature $-1$, so that $h=n-1$, with
$f=f_{\frac{1}{h}\pm\frac{1}{n}}$ where $n\ra+\infty$, Theorem
\ref{theo:introapproxpointconstcurv} of the introduction follows, in
the standard way one deduces a logarithm law-type theorem from a
Khintchine-type theorem.


\noindent {\small 
\begin{tabular}{l} 
University of Georgia\\
Department of Mathematics\\
Athens, GA 30602, USA\\
{\it e-mail: saarh@math.uga.edu} 
\end{tabular} \hfill and\hfill
\begin{tabular}{l} 
Ben Gurion University\\
Department of Mathematics\\
Bear Sheva, Israel\\
{\it e-mail: saarh@cs.bgn.ac.il} 
\end{tabular} 
\\ 
 \mbox{} 
\\ 
 \mbox{} 
\\ 
\begin{tabular}{l}  
D\'epartement de Math\'ematique et Applications, UMR 8553 CNRS\\ 
\'Ecole Normale Sup\'erieure, 45 rue d'Ulm\\ 
75230 PARIS Cedex 05, FRANCE\\ 
{\it e-mail: Frederic.Paulin@ens.fr} 
\end{tabular} 
}

\end{document}